\documentclass[preprint]{imsart}

\RequirePackage[OT1]{fontenc}
\RequirePackage{amsthm,amsmath,amssymb,dsfont,algorithm,algorithmic}
\RequirePackage[numbers]{natbib}
\RequirePackage[colorlinks,citecolor=blue,urlcolor=blue]{hyperref}

\arxiv{arXiv:0000.0000}

\startlocaldefs
\numberwithin{equation}{section}
\theoremstyle{plain}
\newtheorem{thm}{Theorem}[section]
\newtheorem{asm}{Assumption}[section]
\newtheorem{lemma}[thm]{Lemma}
\newtheorem{cor}[thm]{Corollary}
\newtheorem{prop}[thm]{Proposition}
\newtheorem{dfn}{Definition}[section]
\newtheorem{rmk}{Remark}[section]
\DeclareMathOperator*{\argmin}{arg\,min}

\endlocaldefs

\allowdisplaybreaks

\begin{document}

\begin{frontmatter}
\title{Concentration of tempered posteriors and of their variational approximations}
\runtitle{Concentration of variational approximations}

\begin{aug}
\author{\fnms{Pierre} \snm{Alquier}\thanksref{t1}\ead[label=e1]{pierre.alquier@ensae.fr}}
\and
\author{\fnms{James} \snm{Ridgway}\ead[label=e2]{james.lp.ridgway@gmail.com}}

\thankstext{t1}{This author gratefully acknowledges financial support from the research programme {\it New Challenges for New Data} from LCL and GENES, hosted by the {\it Fondation du Risque} and from Labex ECODEC (ANR-11-LABEX-0047).}

\runauthor{P. Alquier and J. Ridgway}

\affiliation{CREST, ENSAE\thanksmark{m1} and Capital Fund Management\thanksmark{m2}}

\address{Pierre Alquier\\
CREST, ENSAE\\
5 avenue Henry le Ch\^atelier
\\
91120 Palaiseau
\\
FRANCE
\\
\printead{e1}}

\address{James Ridgway\\
Capital Fund Management
\\
23 rue de l'Universit\'e
\\
75007 Paris
\\
FRANCE
\\
\printead{e2}}
\end{aug}

\begin{abstract}
While Bayesian methods are extremely popular in statistics and machine learning, their application to massive datasets is often challenging, when possible at all. The classical MCMC algorithms are prohibitively slow when both the model dimension and the sample size are large. Variational Bayesian methods aim at approximating the posterior by a distribution in a tractable family $\mathcal{F}$. Thus, MCMC are replaced by an optimization algorithm which is orders of magnitude faster. VB methods have been applied in such computationally demanding applications as collaborative filtering, image and video processing, or NLP to name a few. However, despite nice results in practice, the theoretical properties of these approximations are not known. We propose a general oracle inequality that relates the quality of the VB approximation to the prior $\pi$ and to the structure of $\mathcal{F}$. We provide a simple condition that allows to derive rates of convergence from this oracle inequality. We apply our theory to various examples. First, we show that for parametric models with log-Lipschitz likelihood, Gaussian VB leads to efficient algorithms and consistent estimators. We then study a high-dimensional example: matrix completion, and a nonparametric example: density estimation.
\end{abstract}



\end{frontmatter}

\section{Introduction}

\subsection{Motivation}

In many applications of Bayesian statistics, the posterior is not tractable. Markov Chain Monte Carlo algorithms (MCMC) were developed to allow the statistician to sample from the posterior distribution even in situations where a closed-form expression is not available. MCMC methods were successfully used in many applications, and are still one of the most valuable tools in the statistician's toolbox. However, many modern applications of statistics and machine learning involve such massive datasets that sampling schemes such as MCMC have become impractical. In order to allow the use of Bayesian approaches with these datasets, it is actually much faster  to compute variational approximations of the posterior by using optimization algorithms. Variational Bayes (VB) has indeed become a corner stone algorithm for fast Bayesian inference. 

VB has been applied to many challenging problems: matrix completion for collaborative filtering~\cite{Lim2007}, NLP on massive datasets~\cite{hoffman2013stochastic}, video processing~\cite{li2009patch}, classification with Gaussian processes~ \cite{gibbs2000variational}, among others. Chapter 10 in~\cite{Bishop2006} is a good introduction to VB and~\cite{blei2017variational} provides an exhaustive survey.

Despite its practical success very little attention has been put towards theoretical guaranties for VB. Asymptotic results in exponential models were provided in \cite{wang2004convergence}. More recently, \cite{wang2017consistency} proposed a very nice asymptotic study of approximations in parametric models. The main problem with these results is that by nature they cannot be applied to high-dimensional or nonparametric models, or to model selection. In the machine learning community,~\cite{alquier2016properties} also studied VB approximations. In a distribution-free setting, there is actually no likelihood, but a pseudo-likelihood can be defined through a suitable loss function and thus it is possible to define a pseudo-posterior. Thanks to PAC-Bayesian inequalities from~\cite{catoni2004statistical,MR2483528},~\cite{alquier2016properties} derived rates of convergence for VB approximation of this pseudo-posterior. However, the tools used in~\cite{alquier2016properties} are valid for bounded loss functions, so there is no direct way to adapt this method to study VB approximations when the log-likelihood is unbounded.

In this paper, we propose a general way to derive concentration rates for approximations of fractional posteriors. Concentration rates are the most natural way to assess ``frequentist guarantees for Bayesian estimators'': the objective is to prove that the posterior is asymptotically highly concentrated around the true value of the parameter. This approach is now very well understood, we refer the reader to the milestone paper~\cite{ghosal2000convergence}, an account of recent advances can be found in in~\cite{rousseau2016frequentist,ghoshal2017book}. Recently,~\cite{bhattacharya2016bayesian} studied the situation where the likelihood $L(\theta)$ is replaced by $L^{\alpha}(\theta)$ for $0<\alpha<1$, leading to what is usually called a {\it fractional} or {\it tempered} posterior. They proved that concentration of the fractional posterior requires actually fewer hypothesis than concentration of the (true) posterior. Extending the technique of~\cite{bhattacharya2016bayesian}, we analyze the concentration of VB approximations of (fractional) posteriors. Especially, we derive a condition for the VB approximation to concentrate at the same rate as the fractional posterior.

\subsection{Definitions and notations}

We observe a collection of $n$ i.i.d. random variables $(X_1,\cdots,X_n)=X_1^n$ in a measured sample space $(\mathbb{X},\mathcal{X},\mathbb{P})$. Let $\left\{P_\theta,\theta\in \Theta\right\}$ be a statistical model (a collection of probability distributions). The objective here is to estimate the distribution of the $X_i$'s. Most results will be stated under the assumption that the model is well specified, {\it i.e.} there exists $\theta_0\in\Theta$ such that $\mathbb{P}\equiv P_{\theta_0}^{\otimes n}$. However, we will also provide results in the case $\mathbb{P}\equiv (P^*)^{\otimes n}$ where $P^*$ does not belong to the model. Let us first assume that $\mathbb{P}\equiv P_{\theta_0}^{\otimes n}$ (we will explicitly mention when this will no longer be the case).

Assume that $Q$ is a dominating measure for this family of distributions, and put $p_{\theta} = \frac{{\rm d}P_{\theta}}{{\rm d} Q}(\theta)$. Let $\mathcal{M}_{1}^+(E)$ be the set of all probability distributions on a measurable space $(E,\mathcal{E})$. Assume $\Theta$ is equipped with some $\sigma$-algebra $\mathcal{T}$. Let $\pi\in \mathcal{M}_{1}^+(\Theta)$ denote the prior. The likelihood and the negative log-likelihood ratio will be denoted respectively\footnote{In order to manipulate these quantities we need to assume that $(X_1^n,\theta)\mapsto r_n(\theta,\theta_0)$ is measurable for the product $\sigma$-field $\mathcal{X}\otimes\mathcal{T}$. This imposes some regularity on $\Theta\times\mathbb{X}$ that will be implicitly assumed in the rest of the paper.} by
\begin{equation*}
 \forall (\theta,\theta')\in\Theta^2\text{, }
 L_n(\theta) = \prod_{i=1}^n p_{\theta}(X_i) \text{ and }
 r_n(\theta,\theta') = \sum_{i=1}^n \log \frac{p_{\theta'}(X_i)}{p_{\theta}(X_i)}.
\end{equation*}
\begin{dfn} Let $\alpha\in(0,1)$. Let $P$ and $R$ be two probability measures. Let $\mu$ be any measure such that $P\ll \mu$ and $R\ll \mu$, for example $\mu=P+R$. The $\alpha$-R\'enyi divergence and the Kullback-Leibler (KL) divergence between two probability distributions $P$ and $R$ are respectively defined by
\begin{align*}
 D_{\alpha}(P,R) & =
 \frac{1}{\alpha-1} \log \int \left(\frac{{\rm d}P}{{\rm d}\mu}\right)^\alpha \left(\frac{{\rm d}R}{{\rm d}\mu}\right)^{1-\alpha} {\rm d}\mu  \text{,}
       \\
 \mathcal{K}(P,R) & = 
 \int \log \left(\frac{{\rm d}P}{{\rm d}R} \right){\rm d}P \text{ if } P \ll R
\text{, }
 + \infty \text{ otherwise}.
\end{align*}
\end{dfn}

\begin{rmk}
\label{properties-renyi}
We remind the reader of a few properties proven in \cite{van2014renyi}. First, it is obvious that $D_{\alpha}(P,R)$ does actually not depend on the choice of the reference measure $\mu$. This is sometimes made explicit by the (informal) statement $ D_{\alpha}(P,R) = (1/(\alpha-1))\log \int ({\rm d}P)^{\alpha} ({\rm d} R)^{1-\alpha}$. The measures $P$ and $R$ are mutually singular if and only if $D_{\alpha}(P,R)=(\frac{1}{\alpha-1})\log(0) = +\infty$.

We have $ \lim_{\alpha\rightarrow 1} D_{\alpha}(P,R) = \mathcal{K}(P,R)$
which gives ground to the notation $D_{1}(P,R)=\mathcal{K}(P,R)$. For $\alpha\in(0,1]$,
$ (\alpha/2) d^{2}_{TV}(P,R) \leq D_{\alpha}(P,R) $, $d_{TV}$ being the total variation distance -- for $\alpha=1$ this is Pinsker's inequality. The map $\alpha\mapsto D_{\alpha}(P,R)$ is nondecreasing. Also, the authors of~\cite{bhattacharya2016bayesian} note that the $\alpha$-R\'enyi divergences are all equivalent for $0<\alpha<1$, through the formula $ \frac{\alpha}{\beta}\frac{1-\beta}{1-\alpha} D_{\beta} \leq D_\alpha \leq D_\beta $ for $\alpha \leq \beta $. Additivity holds:
$
 D_{\alpha}(P_1 \otimes P_2,R_1 \otimes R_2) = D_{\alpha}(P_1,R_1)+D_{\alpha}(P_2,R_2) $, thus $ D_{\alpha}(P^{\otimes n},R^{\otimes n}) = n D_{\alpha}(P,R) $; $D_{1/2}(P,R) \geq 2[1-\exp(-(1/2)D_{1/2}(P,R))] = H^2(P,R)$ the squared Hellinger distance.
\end{rmk}
The fractional posterior, that will be our \textit{ideal} estimator, is given by
\begin{equation*}
\pi_{n,\alpha}({\rm d}\theta|X_1^n):=\frac{{\rm e}^{-\alpha r_n(\theta,\theta_0)} \pi({\rm d}\theta)}{\int{\rm e}^{-\alpha r_n(\theta,\theta_0)} \pi({\rm d}\theta)} \propto L_n^{\alpha}(\theta) \pi({\rm d}\theta),
\end{equation*}
using the notation of \cite{bhattacharya2016bayesian}. The variational approximation $\tilde{\pi}_{n,\alpha}({\rm d}\theta|X_1^n)$ of $\pi_{n,\alpha}({\rm d}\theta|X_1^n)$ is defined as the projection in KL divergence onto a predefined family of distributions $\mathcal{F}$.  
\begin{dfn}[Variational Bayes approximation]
Let $\mathcal{F}\subset \mathcal{M}_{1}^+(\Theta)$,
\begin{equation*}
 \tilde{\pi}_{n,\alpha}(\cdot|X_1^n) = \underset{\rho \in \mathcal{F}}{\arg\min} \, \mathcal{K}(\rho,\pi_{n,\alpha}(\cdot|X_1^n)).
\end{equation*}
\end{dfn}
In Section~\ref{sec:main-res} we state general theorems on the concentration of $\tilde{\pi}_{n,\alpha}(\cdot|X_1^n)$, {\it e.g.} Theorem~\ref{thm-concentration-vb}. One of the key assumptions is that the prior gives enough mass to neighborhoods of the true parameter, a condition also required to prove the concentration of the posterior \cite{ghosal2000convergence,rousseau2016frequentist,bhattacharya2016bayesian}. Here, an additional, but completely natural assumption is required: $\mathcal{F}$ must actually contain distributions concentrated around the true parameter. The choice of $\mathcal{F}$ has thus a strong influence on the quality of the approximation. On one end of the spectrum $\mathcal{F}=\mathcal{M}^1_+(\Theta)$ leads to $\tilde{\pi}_{n,\alpha}=\pi_{n,\alpha}$ and in this situation, our result exactly coincides with the known results on $\pi_{n,\alpha}$. But this is of little interest when $\pi_{n,\alpha}$ is not tractable. On the other end, any family consisting of too few measures will not be rich enough to ensure concentration.

In Sections \ref{sec:gaussian-vb}, \ref{sec:matrix-compl} and \ref{sec:density} we apply our general results in various settings. In Section \ref{sec:gaussian-vb} we study the parametric family of Gaussian approximations
\[
\mathcal{F}^{\Phi}:=\left\lbrace\Phi({\rm d}\theta; m,\Sigma),\quad m\in\mathbb{R}^d,\Sigma\in\mathcal{S}^d_+(\mathbb{R})\right\rbrace 
\]
where $\Phi({\rm d}\theta;m,\Sigma)$ is the $d$ dimensional Gaussian measure with mean $m$ and covariance matrix $\Sigma$, $\mathcal{S}^d_+(\mathbb{R})$ the cone of $d\times d$ symmetric positive definite matrices. We show that other approximations are possible, {\it i.e.} by constraining the variance of the approximation to be a positive diagonal matrix $\Sigma\in\text{Diag}_+^d(\mathbb{R})$. Gaussian approximations have been studied in \cite{Titsias2014,opper2009variational}. We specify those results in the case of a logistic regression in Subsection \ref{sec:logistic}. There, the VB approximation actually turns out to be a convex minimization problem which can be solved by gradient descent or more sophisticated iterative procedures. This is especially attractive as it allows to prove the concentration of the VB approximation obtained after a finite number of steps. In Section \ref{sec:matrix-compl}  we study the case of mean field approximations corresponding to  block-independent distributions
\begin{multline*}
\mathcal{F}^{\text{mf}}:=\bigg\lbrace \rho({\rm d} \theta)=\bigotimes_{i=1}^p\rho_i({\rm d}\theta_i)\in\mathcal{M}^+_1(\Theta),\\
\quad \forall i=1,\cdots,p\quad \rho_i\in\mathcal{M}^+_1(\Theta_i),\quad \Theta=\Theta_1\times\cdots\times\Theta_p\bigg\rbrace,
\end{multline*}
in the context of matrix completion. While the VB approximation leads to feasible approximation algorithms~\cite{Lim2007}, our theorem shows that $\tilde{\pi}_{n,\alpha}$ concentrates at the minimax-optimal rate. In Section \ref{sec:density} we provide a nonparametric example: density estimation. The more important proofs are gathered in Section~\ref{section:proofs}. The supplementary material contains the remaining proofs and additional comments.

\section{Main results}
\label{sec:main-res}

\subsection{A PAC-Bayesian inequality}

We start with a variant of a result of~\cite{bhattacharya2016bayesian}.
\begin{thm}
 \label{thm-bha}
For any $\alpha\in(0,1)$, for any $\varepsilon\in(0,1)$,
 \begin{multline*}
 \mathbb{P}\Biggl( \forall \rho\in \mathcal{M}_{1}^+(\Theta)\text{, }
 \int  D_{\alpha}(P_{\theta},P_{\theta_0}) \rho({\rm d}\theta)\Biggr.
 \\
 \Biggl.\leq \frac{\alpha}{1-\alpha} \int \frac{r_n(\theta,\theta_0) }{n} \rho({\rm d}\theta)
 + \frac{\mathcal{K}(\rho,\pi) + \log\left(\frac{1}{\varepsilon}\right)}{n(1-\alpha)}
 \Biggr) \geq 1-\varepsilon.
 \end{multline*}
\end{thm}
It is tempting to minimize the right-hand side (r.h.s) of the inequality in order to ensure a good estimation.
The minimizer of the r.h.s can actually be explicitly given. In order to do this, let us recall Donsker and Varadhan's variational inequality (Lemma 1.1.3 in~\cite{MR2483528}).
\begin{lemma}
\label{thm-dv}
For any probability $\pi$ on $(\Theta,\mathcal{T})$ and any measurable function
  $h : \Theta \rightarrow \mathbb{R}$ such that $\int{\rm e}^h  \rm{d}\pi < \infty$,
  \begin{equation*}
    \log\int {\rm e}^h \mathrm{d}\pi = \underset{\rho \in \mathcal{M}_{1}^+(\Theta)}{\sup}
    \left[ \int h \mathrm{d}\rho - \mathcal{K}(\rho,\pi) \right],
  \end{equation*}
  with the convention
  $\infty-\infty =-\infty$. Moreover when $h$ is upper-bounded on the
  support of $\pi$ the supremum with respect to $\rho$ in the r.h.s is reached by $\pi_h$ given by $
   \mathrm{d}\pi_h / \mathrm{d}\pi (\theta) = \exp(h(\theta)) / \int \exp(h) \mathrm{d}\pi$.
\end{lemma}
Using Lemma~\ref{thm-dv} with $h(\theta)= - \alpha r_n(\theta,\theta_0)$ and the definition of $\pi_{n,\alpha}$ we obtain
$$
 \pi_{n,\alpha}(\cdot|X_1^n)  = \underset{\rho \in \mathcal{M}^1_+(\Theta)}{\arg\min} \, \left\{
 \alpha \int r_n(\theta,\theta_0) \rho({\rm d}\theta) + \mathcal{K}(\rho,\pi)
 \right\}
$$
so the minimizer of the r.h.s of Theorem~\ref{thm-bha} is actually $\pi_{n,\alpha}({\rm d}\theta|X_1^n)$.

The statement of Theorem~\ref{thm-bha} for $\rho=\pi_{n,\alpha}({\rm d}\theta|X_1^n)$ is Theorem 3.5 in~\cite{bhattacharya2016bayesian}. The proof of Theorem~\ref{thm-bha} requires a straightforward extension, we provide it in Section~\ref{section:proofs} for the sake of completeness. Our extension is crucial though as we will have to use it with $\rho=\tilde{\pi}_{n,\alpha}({\rm d}\theta|X_1^n)$.

\begin{rmk} \label{remark-laplace}
 Theorem~\ref{thm-bha} can be used to study other approximations of the posterior. For example, as suggested by one of the Referees, we can use it to study distributions centered around the {\it maximum a posteriori} (MAP) or the maximum likelihood estimate (MLE). For example, Laplace approximations are Gaussian distributions centered at the MLE. However, there are models $(P_\theta,\theta\in\Theta)$ where the MLE and the MAP are not defined, while the posterior and some variational approximations are consistent. Such an example is provided in the Supplementary Material.
\end{rmk}

\subsection{Concentration of VB approximations}
\label{section:vb-concentration}

We specialize the above results to the variational approximation. Elementary calculations show that
\begin{align*}
 \tilde{\pi}_{n,\alpha}(\cdot|X_1^n) & = \underset{\rho \in \mathcal{F}}{\arg\min} \, \left\{
 \alpha \int r_n(\theta,\theta_0) \rho({\rm d}\theta) + \mathcal{K}(\rho,\pi)
 \right\}
 \\
 & = \underset{\rho \in \mathcal{F}}{\arg\min} \, \left\{
 - \alpha \int \sum_{i=1}^n \log p_{\theta}(X_i) \rho({\rm d}\theta) + \mathcal{K}(\rho,\pi)
 \right\}.
\end{align*}
As a consequence, we obtain the following corollary of Theorem~\ref{thm-bha}.
\begin{cor}
\label{thm-cor-bha}
 For any $\alpha\in(0,1)$ and $\varepsilon\in(0,1)$, with probability at least $1-\varepsilon$,
 \begin{align*}
 & \int D_{\alpha}(P_{\theta},P_{\theta_0}) \tilde{\pi}_{n,\alpha}({\rm d}\theta|X_1^n)
  \\
 & \quad \leq \inf_{\rho\in\mathcal{F}} \left\{ \frac{\alpha}{1-\alpha} \int \frac{r_n(\theta,\theta_0) }{n} \rho({\rm d}\theta)
 + \frac{\mathcal{K}(\rho,\pi) + \log\left(\frac{1}{\varepsilon}\right)}{n(1-\alpha)}
 \right\}.
 \end{align*}
\end{cor}
Obviously, when $\mathcal{F}=\mathcal{M}_{1}^+(\Theta)$, we have $\tilde{\pi}_{n,\alpha}({\rm d}\theta|X_1^n)=\pi_{n,\alpha}({\rm d}\theta|X_1^n)$, so we recover as a special case an upper bound on the risk of the tempered posterior. We are now in position to state our main result.
\begin{thm}
\label{thm-concentration-vb}
Fix $\mathcal{F}\subset \mathcal{M}_{1}^+(\Theta)$.
 Assume that a sequence $\varepsilon_n>0$ is such that there is a distribution $\rho_n\in\mathcal{F}$ such that
  \begin{equation}
  \label{equa-assumption-1}
  \int \mathcal{K}(P_{\theta_0},P_{\theta}) \rho_n({\rm d}\theta) \leq \varepsilon_n \text{, }
  \int \mathbb{E} \left[ \log^2 \left( \frac{p_{\theta}(X_i)}{p_{\theta_0}(X_i)} \right) \right]  \rho_n({\rm d}\theta) \leq \varepsilon_n
 \end{equation}
 and
 \begin{equation}
   \label{equa-assumption-2}
  \mathcal{K}(\rho_n,\pi) \leq n \varepsilon_n.
 \end{equation}
 Then, for any $\alpha\in(0,1)$, for any $(\varepsilon,\eta)\in(0,1)^2$,
  \begin{equation*}
  \mathbb{P}\left[
 \int D_{\alpha}(P_{\theta},P_{\theta_0}) \tilde{\pi}_{n,\alpha}({\rm d}\theta|X_1^n)  \leq  \frac{(\alpha+1) \varepsilon_n + \alpha \sqrt{\frac{\varepsilon_n}{n\eta}} + \frac{\log\left(\frac{1}{\varepsilon}\right)}{n}}{1-\alpha}\right]\geq 1-\varepsilon-\eta.
 \end{equation*}
\end{thm}
This theorem is a consequence of Corollary~\ref{thm-cor-bha}, its proof is provided in Section~\ref{section:proofs}. Let us now discuss the main consequences of this theorem.

Note that the assumption involving a distribution $\rho_n$ is not standard. This requires some explanations. Consider first the case $\mathcal{F}=\mathcal{M}_{1}^+(\Theta)$. Define $B(r)$, for $r>0$, as
 \begin{equation*}
  B(r) = \left\{ \theta\in\Theta: \mathcal{K}(P_{\theta_0},P_{\theta}) \leq r, \mathbb{E} \left[ \log^2 \left( \frac{p_{\theta}(X_i)}{p_{\theta_0}(X_i)} \right) \right]  \leq r \right\}.
 \end{equation*}
 Then the choice $\rho_n=\pi_{|B(\varepsilon_n)}$, {\it i.e.} $\pi$ restricted to $B(\varepsilon_n)$, ensures immediately~\eqref{equa-assumption-1}, and~\eqref{equa-assumption-2} can be rewritten
 \begin{equation*}
  - \log \pi(B(\varepsilon_n)) \leq n \varepsilon_n.
 \end{equation*}
 This assumption is standard to study concentration of the posterior, see Theorem 2.1 page 503 in~\cite{ghosal2000convergence} or Subsection 3.2 in~\cite{rousseau2016frequentist}. Our message is that in the studies of concentration of the posterior, the choice $\rho_n=\pi_{|B(\varepsilon_n)}$ was hidden. Other choices might lead to easier calculations in some situations. More importantly, in the relevant case $\mathcal{F}\subsetneq \mathcal{M}_{1}^+(\Theta)$, $\pi_{|B(\varepsilon_n)}\notin\mathcal{F}$ in general. Thus $- \log \pi(B(\varepsilon_n)) \leq n \varepsilon_n$ is no longer sufficient, and~\eqref{equa-assumption-1} and~\eqref{equa-assumption-2} are natural extensions of this assumption to study VB. They provide an explicit condition on the family $\mathcal{F}$ in order to ensure concentration of the approximation.

Choosing $\eta=\frac{1}{n \varepsilon_n} $ and $\varepsilon = \exp(-n \varepsilon_n)$ we obtain a more readable concentration result. It shows that, as soon as $(1/n) \ll \varepsilon_n \ll 1$, the sequence $\varepsilon_n$ gives a concentration rate for VB.
\begin{cor}
\label{cor-concentration}
Under the same assumptions as in Theorem~\ref{thm-concentration-vb},
  \begin{align*}
  \mathbb{P}\left[
 \int D_{\alpha}(P_{\theta},P_{\theta_0}) \tilde{\pi}_{n,\alpha}({\rm d}\theta|X_1^n)  \leq \frac{2(\alpha+1)}{1-\alpha} \varepsilon_n\right] & \geq 1-\frac{1}{n\varepsilon_n}-\exp(-n\varepsilon_n)
 \\
 & \geq 1-\frac{2}{n\varepsilon_n},
 \end{align*}
\end{cor}

\begin{rmk} 
 As a special case, when $\alpha=1/2$, the theorem leads to a concentration result in terms of the more classical Hellinger distance
   \begin{equation*}
  \mathbb{P}\left[
 \int H^2(P_{\theta},P_{\theta_0}) \tilde{\pi}_{n,1/2}({\rm d}\theta|X_1^n)  \leq 6 \varepsilon_n \right]\geq 1-\frac{2}{n\varepsilon_n}.
 \end{equation*}
 Also, with a general $\alpha \in(0,1)$, from the properties recalled in Remark~\ref{properties-renyi}, we have, for $0<\beta\leq \alpha$,
$$
  \mathbb{P}\left[
 \int D_{\beta}(P_{\theta},P_{\theta_0}) \tilde{\pi}_{n,\alpha}({\rm d}\theta|X_1^n)  \leq \frac{2(\alpha+1)}{1-\alpha} \varepsilon_n\right]
 \geq 1-\frac{2}{n\varepsilon_n},
$$
and for $\alpha \leq \beta < 1$,
 $$
  \mathbb{P}\left[
 \int D_{\beta}(P_{\theta},P_{\theta_0}) \tilde{\pi}_{n,\alpha}({\rm d}\theta|X_1^n)  \leq \frac{2\beta (\alpha+1)}{\alpha(1-\beta)} \varepsilon_n\right]
 \geq 1-\frac{2}{n\varepsilon_n}.
$$
\end{rmk}
 
\subsection{A simpler result in expectation}

It is possible to simplify the assumptions at the price of stating a result in expectation instead of concentration.

\begin{thm}
\label{thm-concentration-vb:expect}
Fix $\mathcal{F}\subset \mathcal{M}_{1}^+(\Theta)$. Then
\begin{multline*}
\mathbb{E} \left[ \int D_{\alpha}(P_{\theta},P_{\theta_0}) \tilde{\pi}_{n,\alpha}({\rm d}\theta|X_1^n) \right]
\\
\leq \inf_{\rho\in\mathcal{F}} \left\{ \frac{\alpha}{1-\alpha} \int \mathcal{K}(P_{\theta_0},P_{\theta}) \rho({\rm d}\theta)
 + \frac{\mathcal{K}(\rho,\pi) }{n(1-\alpha)}
 \right\}.
\end{multline*}
Assume that $\varepsilon_n>0$ is such that there is distribution $\rho_n\in\mathcal{F}$ such that
$$
  \int \mathcal{K}(P_{\theta_0},P_{\theta}) \rho_n({\rm d}\theta)  \leq\varepsilon_n
  \text{, and }\mathcal{K}(\rho_n,\pi)  \leq n\varepsilon_n.
$$
 Then, for any $\alpha\in(0,1)$,
\begin{equation*}
\mathbb{E} \left[ \int D_{\alpha}(P_{\theta},P_{\theta_0}) \tilde{\pi}_{n,\alpha}({\rm d}\theta|X_1^n) \right]
\leq \frac{1+\alpha}{1-\alpha}\varepsilon_n.
\end{equation*}
\end{thm}
 
\subsection{Extension of the result in expectation to the misspecified case}

In this section we do not assume any longer that the true distribution is in $\{P_{\theta},\theta\in\Theta\}$. In order not to change all the notations we define an extended parameter set $\Theta \cup \{\theta_0\}$ where $\theta_0\notin\Theta$ and define $P_{\theta_0}$ as the true distribution. Theorem~\ref{thm-concentration-vb:expect} can be applied to this setting, and we obtain:
\begin{multline*}
\mathbb{E} \left[ \int D_{\alpha}(P_{\theta},P_{\theta_0}) \tilde{\pi}_{n,\alpha}({\rm d}\theta|X_1^n) \right]
\\
\leq \inf_{\rho\in\mathcal{F}} \left\{ \frac{\alpha}{1-\alpha} \int  \mathcal{K}(P_{\theta_0},P_{\theta})   \rho({\rm d}\theta)
 + \frac{\mathcal{K}(\rho,\pi) }{n(1-\alpha)}
 \right\}.
\end{multline*}
Now, rewriting, for $\theta^*\in\Theta$,
$$
 \mathcal{K}(P_{\theta_0},P_{\theta}) = \mathcal{K}(P_{\theta_0},P_{\theta^*}) + \mathbb{E}\left[\log\frac{ p_{\theta^*}(X_i)}{p_{\theta}(X_i)}\right],
$$
we obtain the following result.
\begin{thm}
\label{thm-concentration-vb:expect:misspecified}
Assume that, for $\theta^*=\arg\min_{\theta\in\Theta} \mathcal{K}(P_{\theta_0},P_\theta)$, there is $\varepsilon_n>0$ and $\rho_n\in\mathcal{F}$ with
$$
  \int \mathbb{E}_{\theta_0}\left[\log\frac{ p_{\theta^*}(X_i)}{ p_{\theta}(X_i)}\right]\rho_n({\rm d}\theta)   \leq \varepsilon_n
  \text{, and } \mathcal{K}(\rho_n,\pi) \leq n \varepsilon_n,
$$
 then, for any $\alpha\in(0,1)$,
\begin{equation*}
\mathbb{E} \left[ \int D_{\alpha}(P_{\theta},P_{\theta_0}) \tilde{\pi}_{n,\alpha}({\rm d}\theta|X_1^n)\right]
\leq \frac{\alpha}{1-\alpha} \min_{\theta\in\Theta} \mathcal{K}(P_{\theta_0},P_{\theta})+ \frac{1+\alpha}{1-\alpha}\varepsilon_n.
\end{equation*}
\end{thm}
In the well-specified case, $\theta^*=\theta_0$ and we recover Theorem~\ref{thm-concentration-vb:expect}. Otherwise, this result takes the form of an oracle inequality. It is not a sharp oracle inequality as that the risk measure used in the l.h.s and the r.h.s are not the same, but remains informative when $\mathcal{K}(P_{\theta_0},P_{\theta^*})$ is small. For example, in Section~\ref{sec:density} below, we provide a nonparametric example where $\mathcal{K}(P_{\theta_0},P_{\theta^*})$ and Theorem~\ref{thm-concentration-vb:expect:misspecified} leads to the minimax rate of convergence.

\section{Gaussian variational Bayes}
\label{sec:gaussian-vb}

In this section we consider $\Theta\subset\mathbb{R}^d$ and the class of Gaussian approximations 
$$\mathcal{F}^\Phi:=\left\lbrace \Phi(d\theta;m,\Sigma),\quad m\in\mathbb{R}^d,\Sigma\in\mathcal{S}^d_+(\mathbb{R})\right\rbrace , $$
thus the algorithm will consist in projecting onto the set of Gaussian distributions.
Depending on the hypotheses made on the covariance matrix we can build different approximations. For instance define:
\begin{align*}
  \mathcal{F}^\Phi_{diag}&:=\left\lbrace \Phi(d\theta;m,\Sigma),\quad m\in\mathbb{R}^d,\Sigma\in\text{Diag}_+^d(\mathbb{R})\right\rbrace\\
  \mathcal{F}^\Phi_{id}&:=\left\lbrace \Phi(d\theta;m,\sigma^2I_d),\quad m\in\mathbb{R}^d,\sigma^2\in\mathbb{R_{+\star}}\right\rbrace.
\end{align*}
We have by definition $\mathcal{F}^\Phi_{id}\subseteq\mathcal{F}^\Phi_{diag}\subseteq\mathcal{F}^\Phi$. 

The remarkable fact of Gaussian VB is that it allows to recast integration as a finite dimension optimization problem. 
The choice of a specific Gaussian is a trade off between accuracy and computational complexity. We will show in the following that, 
under some assumption on the likelihood, the integrated $\alpha$-R\'enyi  divergence is convergent for most of the approximations.

To simplify the exposition of the results we will restrict our study to the case of Gaussian priors: $\pi=\mathcal{N}(0,\vartheta^2 I_p)$. One can readily see that in
Theorem \ref{thm-concentration-vb} the prior appears only in the condition $\frac1n\mathcal{K}(\rho,\pi)\leq \varepsilon_n$,
 many other distribution could be used, providing different rates.

In the rest of the section we assume that the density is log Lipschitz.  
\begin{asm}
\label{lipschitz-den}
There is a measurable real function $M(\cdot)$ such that
\[
\left\vert\log p_\theta(X_1)-\log p_{\theta^\prime}(X_1) \right\vert \leq M(X_1)\left\Vert\theta-\theta^\prime\right\Vert_2
\]
Furthermore we assume that $\mathbb{E}M(X_1)=:B_1,\quad\mathbb{E}M^2(X_1)=:B_2<\infty$.
\end{asm}
An example is logistic regression, see Subsection~\ref{sec:logistic} below.
  \begin{thm}
\label{thm-lipschitz-dens}
  Let the approximation family be $\mathcal{F}$ with $\mathcal{F}^\Phi_{id}\subset \mathcal{F}$ as defined above
    and that the model satisfies Assumption \ref{lipschitz-den}. We put
    $$
\varepsilon_n=\frac{B_1}n\vee \frac{B_2}{n^2} \vee \left\lbrace\frac dn \left[\frac12
  \log\left(\vartheta^2 n^2 \sqrt{d} \right)
        + \frac{1}{n\vartheta^2}
  \right]
  + \frac{\|\theta_0\|^2}{n\vartheta^2} - \frac{d}{2n}\right\rbrace
$$
Then for any
    $\alpha\in(0,1)$, for any $\eta,\epsilon$
    \begin{equation*}
    \mathbb{P}\left[
 \int D_{\alpha}(P_{\theta},P_{\theta_0}) \tilde{\pi}_{n,\alpha}({\rm d}\theta|X_1^n)  \leq  \frac{(\alpha+1) \varepsilon_n + \alpha \sqrt{\frac{\varepsilon_n}{n\eta}} + \frac{\log\left(\frac{1}{\varepsilon}\right)}{n}}{1-\alpha}\right]\geq 1-\varepsilon-\eta.
 \end{equation*}
  \end{thm}
\subsection{Stochastic variational Bayes}
In many cases the model is not conjugate, {\it i.e.} the VB objective does not have a closed-form solution. 
We can however use a full Gaussian approximation and a stochastic gradient descent
on the objective function defined by the KL divergence. This approach has been studied in \cite{Titsias2014}.

We may write our variational bound as the following minimisation problem
$$
\min_{\rho\in\mathcal{F}_\Phi}\int \rho(d\theta)\log \frac{d\rho(\theta)}{d\pi_{n,\alpha}(\theta\vert X^n)}
$$
or after dropping the constants,
\begin{equation}
\label{vb-objective}
\min_{m\in\mathbb{R}^d,\Sigma\in\mathcal{S}^d_+}\left\lbrace-\alpha\int\log p_\theta(Y^n)\Phi(d\theta;m,\Sigma)+\int\log \frac{d\Phi(\theta;m,\Sigma)}{d\pi(\theta)}\Phi(d\theta;m,\Sigma)\right\rbrace.
\end{equation}

In \cite{Titsias2014} the authors suggest using a parametrization of the problem where we replace the optimization over $\Sigma$ by a minimization over the 
matrix $C$ where $CC^t=\Sigma$. To simplify the notations in this section define 
$$
F:x=\left(m,C\right)\in \mathbb{R}^d\times\mathbb{R}^{d\times d}\mapsto\mathbb{E}\left[f(x,\xi)\right]
$$
to be the objective of the minimization problem \eqref{vb-objective}, where $\xi\sim\mathcal{N}(0,I_d)$ and 
\begin{equation}
 \label{eq:dfn_f}
f((m,C),\xi):=-\alpha \log p_{m+C\xi}(Y_1^n)+\log\frac{{\rm d}\Phi_{m,CC^t}}{{\rm d}\pi}(m+C\xi).
\end{equation}
In order to be able to state non-asymptotic results on the stochastic gradient algorithms, we restrict the parameter space to an Euclidean ball, that is \eqref{vb-objective} is transformed into
\[
\min_{x\in\mathbb{B}\cap\mathbb{R}^d\times\mathcal{R}^{d\times d}}\mathbb{E}\left[f(x,\xi)\right],
\]
where $\mathbb{B}=\left\lbrace x\in\mathbb{R}^{d^2+d}, \left\Vert x\right\Vert_2\leq B \right\rbrace$ for some $B>0$. We will then let $\mathcal{P}_\mathbb{B}$ denote the orthogonal projection onto $\mathbb{B}$. In addition we can define the corresponding 
family of Gaussian distributions
$$\mathcal{F}_B^\Phi=\left\lbrace \Phi(d\theta;m,CC^t),\quad (m,C)\in \mathbb{B}\cap\mathbb{R}^d\times\mathbb{R}^{d\times d}\right\rbrace .$$

The objective can now be replaced by a Monte Carlo estimate and we can use stochastic gradient descent as described in Algorithm \ref{algo-SVB}.

\begin{algorithm}[h!]
\caption{Stochastic Variational Bayes
\label{algo-SVB}}
\begin{description}
\item[Input:] $x_0$, $X_1^n$, $\gamma_T$ 
\item[For $t\in\left\lbrace1,\cdots,T\right\rbrace$],
\begin{description}
\item[a.] Sample $\xi_t\sim \mathcal{N}(0,I_d)$
\item[b.] Update $x_t\leftarrow \mathcal{P}_\mathbb{B}\left(x_{t-1}-\gamma_T\nabla f(x_{t-1},\xi_t)\right)$
\end{description}
\item[End For].
\item[Output:] $\bar{x}_T=\frac1T\sum_{t=1}^T x_t$
\end{description}
\end{algorithm}

\begin{asm}
\label{asm:conv-lip}
Assume that $f$, as defined in~\eqref{eq:dfn_f}, is convex in its first component $x$ and that it has $L$-Lipschitz gradients.
\end{asm}

Define $\tilde{\pi}_{n,\alpha}^k({\rm d}\theta\vert X_1^n)$ to be the $k$-th iterate of Algorithm \ref{algo-SVB},
the Gaussian distribution with parameters $\bar{x}_k=(\bar{m}_k,\bar{C}_k)$.

\begin{thm} 
\label{thm-vb-convex}
Let Assumptions \ref{asm:conv-lip} and  \ref{lipschitz-den} be verified, and define $\varepsilon_n$ as in Theorem~\ref{thm-lipschitz-dens}. Let $B$ be such that $B>\Vert\theta_0\Vert_2+ 1/n\sqrt{d}$. Then for $\tilde{\pi}^ T_{n,\alpha}({\rm d}\theta|X_1^n)$ obtained by 
Algorithm \ref{algo-SVB} with $\gamma_T=\frac{B}{L\sqrt{2T}}$, we get 
$$
\mathbb{E}  \left[  \int D_{\alpha}(P_{\theta},P_{\theta_0})\tilde{\pi}^T_{n,\alpha}({\rm d}\theta|X_1^n) \right]
\leq \frac1{n(1-\alpha)}\sqrt{\frac{2BL}{T}}+\frac{1+\alpha}{1-\alpha}\varepsilon_n.
$$
\end{thm}
\begin{rmk}
On most examples the gradient is a sum of at least $n$ components. If each term is Lipschitz with constant $L_i$, an estimate of the constant will be $L\leq n\max_i L_i$. The additional term of the bound is therefore of the order $(2B\max L_i/(nT))^{1/2}/ (1-\alpha)$, hence a good choice is $T=O(\sqrt{n})$ to mitigate the impact of the numerical approximation on the rate.
\end{rmk}

\subsection{Example: logistic regression}
\label{sec:logistic}
We consider the case of a binary regression model. Although estimation of parameters is relatively simple 
for small datasets \cite{chopin2017}, it remains challenging when the size of the dictionary is large. Furthermore usual deterministic methods do not come with theoretical guarantees as would a gradient descent algorithm for maximum likelihood. Note that the logistic regression is not conjugate in the sense that we cannot find an iterative scheme based on a mean field approximation, as will be done for the matrix completion example in Section~\ref{sec:matrix-compl}.

Let $X_i=(Y_i,Z_i)\in \lbrace -1,1\rbrace\times\mathbb{R}^d$ be such that
$$
\mathbb{P}\lbrace Y=y \vert Z=z,\theta\rbrace=\frac{e^{yz^t \theta}}{1+e^{yz^t\theta}},
$$
 We will consider the case of estimation with a Gaussian prior $\pi(d\theta)=\Phi(d\theta;0,\vartheta I_d)$ on $(\mathbb{R}^d,\mathcal{B}(\mathbb{R}^d))$ (other cases are easily incorporated in the theory).

 We will prove results in the case of random design where we suppose that the distribution of $Z_1^n$ does not depend on the parameter.

  \begin{cor}
\label{cor:logit-cor}
  Let the family of approximation be any $\mathcal{F}$ with $\mathcal{F}^\Phi_{id}\subset \mathcal{F}$ as defined above
    and assume that $K_1:=2\mathbb{E}\left\Vert X_1\right\Vert, K_2:=4\mathbb{E}\left\Vert X_1\right\Vert^2<\infty$. Put
$$
\varepsilon_n=\frac{K_1}n\vee \frac{K_2}{n^2} \vee \left\lbrace\frac dn \left[
  \frac12 \log\left(\vartheta^2 n^2 \sqrt{d}\right)
        + \frac{1}{n\vartheta^2}
  \right]
  + \frac{\|\theta_0\|^2}{n\vartheta^2} - \frac{d}{2n}\right\rbrace
$$
    then for any
    $\alpha\in(0,1)$, for any $\eta,\epsilon$
    \begin{equation*}
    \mathbb{P}\left[
 \int D_{\alpha}(P_{\theta},P_{\theta_0}) \tilde{\pi}_{n,\alpha}({\rm d}\theta|X_1^n)  \leq  \frac{(\alpha+1) \varepsilon_n + \alpha \sqrt{\frac{\varepsilon_n}{n\eta}} + \frac{\log\left(\frac{1}{\varepsilon}\right)}{n}}{1-\alpha}\right]\geq 1-\varepsilon-\eta.
 \end{equation*}
  \end{cor}

To apply Theorem \ref{thm-vb-convex} we need to add some constraint on the covariance matrix. The optimization will be written 
over $\mathbb{B}_\psi:=\mathbb{B}\cap \lbrace C\in \mathbb{R}^{d\times d}, CC^t\succeq\psi I_{d\times d} \rbrace$ (this is done only to ensure that $\log\vert\Sigma\vert$ has Lipschitz gradients). 
\begin{cor} 
\label{cor:optim-logit}
 Let the family of approximation be any $\mathcal{F}$ with $\mathcal{F}^\Phi_{id}\subset \mathcal{F}$
and assume that $K_1:=2\mathbb{E}\left\Vert X_1\right\Vert, K_2:=4\mathbb{E}\left\Vert X_1\right\Vert^2<\infty$,
let $B$ be such that $B>\Vert\theta_0\Vert_2+ \frac1{n\sqrt{d}}$
 then for $\pi^T_{n,\alpha}({\rm d}\theta|X_1^n)$ obtained by 
Algorithm \ref{algo-SVB} with $\gamma_T=\frac{B}{L\sqrt{2T}}$ and where $\mathbb{B}$ is replaced by $\mathbb{B}_\psi$ for any $\psi\leq\frac1{n\sqrt{d}}$, we get 
\begin{multline*}
\mathbb{E}  \left[  \int D_{\alpha}(P_{\theta},P_{\theta_0})\tilde{\pi}^T_{n,\alpha}({\rm d}\theta|X_1^n) \right]
\leq
\frac1{(1-\alpha)n}\sqrt{\frac{2BL}{T}}
\\
+ \frac{1+\alpha}{1-\alpha}\left(\frac{K_1}n\vee \frac{K_2}{n^2} \vee
\left\lbrace\frac dn \left[
  \frac12\log\left(\vartheta^2 n^2 \sqrt{d}\right)
        + \frac{1}{n\vartheta^2\sqrt{d}}
  \right]
  + \frac{\|\theta_0\|^2}{n\vartheta^2} - \frac{d}{2n}\right\rbrace\right).
\end{multline*}
\end{cor}

Note that the only assumption on the distribution of $X_1$ is that $K_2<\infty$. Still, it is interesting to compute $K_1$ and $K_2$ on some examples. For example, when $X_1$ is uniform on the unit sphere, $K_1 \leq 2$ and $K_2 \leq 4$. When $X_1 \sim \mathcal{N}(0,s^2 I_d) $ then $K_2 = 4 s^2 d$ and $K_1 \leq 2 \sqrt{s^2d}$. In both cases, the terms in $K_1$ and $K_2$ do not deterioriate the parametric rate of convergence $d/n$. Furthermore the Lipschitz constant can be bounded explicitly under additional assumptions on the design matrix ({\it e.g.} bounded singular value) and leads to $L=\mathcal{O}\left(nd+\frac d\psi\right)$. Hence taking $\psi=1/(n\sqrt{d})$ one would get a bound in $\mathcal{O}\left(\sqrt{\frac{d^{3/2}}{nT}}+(d/n)\log{nd}\right)$. We can take $T$ of the order $\frac n{d^{1/4}}$ in order not to deteriorate the rate.
\section{Application to matrix completion}
\label{sec:matrix-compl}

\subsection{Context}

Challenging applications such as collaborative filtering made matrix completion one of the most important machine learning problems in the past few years. Let us describe briefly the model: in this case, our parameter $\theta$ is a matrix $M\in\mathbb{R}^{m \times p}$, with $m,p\geq 1$. For clarity, we will denote by $M_0$ the true matrix $\theta_0$ and use $M$ as a notation for a generic parameter instead of $\theta$. Under $P_M$, the observations are random entries of this matrix with possible noise:
\begin{equation*}
Y_k = M_{i_k,j_k} + \xi_k \text{ for } 1\leq k\leq n
\end{equation*}
where the $(i_k,j_k)$ are i.i.d $\mathcal{U}(\{1,\dots,m\}\times\{1,\dots,p\})$. For the sake of simplicity we will assume that the $\xi_k$ are i.i.d $ \mathcal{N}(0,\sigma^2)$, and that $\sigma^2$ is known, so we only have to estimate $M$. Note that for $\alpha\leq 1$, $
\mathcal{D}_{\alpha}(\mathcal{N}(\mu_1,\sigma^2),\mathcal{N}(\mu_2,\sigma^2)) = \alpha(\mu_1-\mu_2)^2/(2\sigma^2) $, see (10) page 3800 in~\cite{van2014renyi}. Thus, for $0<\alpha<1$,
$$
 \mathcal{D}_{\alpha}(P_M,P_N)  = \frac{1}{\alpha-1} \log \left[ \frac{1}{mp}\sum_{i=1}^m\sum_{j=1}^p \exp\left(\frac{\alpha(\alpha-1)(M_{i,j}-N_{i,j})^2}{2\sigma^2}\right) \right]
$$
which depends only on $\alpha$, $\sigma_2$ and the matrices $M$ and $N$ so we will use the notation $d_{\alpha,\sigma}(M,N) = \mathcal{D}_{\alpha}(P_M,P_N)$. In the case $\alpha=1$,
 $$
 \mathcal{K} (P_M,P_N)  = \frac{1}{mp} \sum_{i=1}^m \sum_{j=1}^p \frac{ (M_{i,j}-N_{i,j})^2}{2\sigma^2}= \frac{ \|M-N\|_F^2}{2\sigma^2 mp}
$$
where $\|\cdot\|_F$ denotes the Frobenius norm. In the noiseless case $\sigma^2=0$,~\cite{candes2010power} proved that it is possible to recover exactly $M$ under the assumption that its rank is small enough. Various extensions to noisy settings, approximately low-rank matrices, or other loss functions can be found in~\cite{candes2010matrix,koltchinskii2011nuclear,klopp2014noisy,alquiercottetlecue}. The main message of these papers is that the minimax rate of convergence is $(m+p){\rm rank}(M)/n$, possibly up to log terms. Bayesian estimators were proposed in~\cite{salakhutdinov2008bayesian,lawrence2009non,zhou2010nonparametric,alquier2014bayesian} using factorized Gaussian priors. Convergence of the posterior mean was proven in~\cite{alquier2015bayesian} for a bounded prior, excluding the Gaussian prior used in practice. Similarly,~\cite{suzuki2015convergence} proves concentration of a truncated version of the posterior. For very large datasets the MCMC algorithm proposed in~\cite{salakhutdinov2008bayesian} is too slow, a VB approximation was proposed in~\cite{Lim2007} with very good results on the Netflix dataset. This approximation was re-used and extended by many authors including~\cite{Paisley,Babacan,alquier2014bayesian,cottetalquier,marsden}. But the consistency of the Bayesian estimator with Gaussian priors and of its variational approximations are opened questions.

First, we will recall the Gaussian prior~\cite{salakhutdinov2008bayesian} and the VB approximation~\cite{Lim2007}. We will then prove the concentration of the VB approximation, and as a consequence the concentration of the tempered posterior.

\subsection{Definition of the prior and of the VB approximation}

Fix $K\in\{1,\dots,m\wedge p\}$. The main idea of factorized priors is that, when ${\rm rank}(M) \leq K$ then we have
\begin{equation*}
 M = UV^t
\end{equation*}
for some matrices $U$ of dimension $p\times K$ and $V$ of dimension $m\times K$. Thus, we can define a prior on $M$ by specifying priors on $U$ and $V$. A usual choice is that the entries $U_{i,k}$ and $V_{j,k}$ are independent $\mathcal{N}(0,\gamma_k)$ and finally $\gamma_k$ is inverse gamma, that is $1/\gamma_k \sim \Gamma(a,b)$. These choices ensure conjugacy: put $\gamma=(\gamma_1,\dots,\gamma_K)$, it is then possible to compute the conditional posteriors of $U|V,\gamma$, of $V|U,\gamma$ and $\gamma|U,V$. This allows to use the Gibbs sampler~\cite{salakhutdinov2008bayesian}. For large datasets,~\cite{Lim2007} proposed mean-field VB with $\mathcal{F}$ given by
\begin{equation*}
 \rho({\rm d}U,{\rm d}V,{\rm d}\gamma) = \bigotimes_{i=1}^{m} \rho_{U_i}({\rm d}U_{i,\cdot})
 \bigotimes_{j=1}^{p} \rho_{V_j}({\rm d}V_{j,\cdot}) \bigotimes_{k=1}^{K} \rho_{\gamma_k}(\gamma_k).
\end{equation*}
The minimization of the VB program is shown in many cited papers, see \cite{alquier2014bayesian} and all the references therein. Shortly: $\rho_{U_i} $ is $\mathcal{N}(\mathbf{m}_{i,\cdot}^t,\mathcal{V}_i)$,
$\rho_{V_j}$ is $\mathcal{N}(\mathbf{n}_{j,\cdot}^t,\mathcal{W}_j)$ and
$\rho_{\gamma_k}$ is $\Gamma(a+(m_1+m_2)/2,\beta_k)$ for some $m\times K$ matrix $\mathbf{m}$
whose rows are denoted by $\mathbf{m}_{i,\cdot}$, some $p\times K$ matrix $\mathbf{n}$
whose rows are denoted by $\mathbf{n}_{j,\cdot}$ and some vector $\beta=(\beta_1,\dots,\beta_K)$.
The parameters are updated iteratively through the formulae
\begin{enumerate}
 \item moments of $U$:
 $$\mathbf{m}_{i,\cdot}^t :=\frac{2\alpha}{n} \mathcal{V}_{i}
\sum_{k:i_k=i} Y_{i_k,j_k} \mathbf{n}_{j_k,\cdot}^t$$
$$
 \mathcal{V}_i^{-1} := \frac{2\alpha}{n} \sum_{k:i_k=i}
  \left[\mathcal{W}_{j_k} + \mathbf{n}_{j_k,\cdot} \mathbf{n}_{j_k,\cdot}^t \right]
    + \left(a+\frac{m_1+m_2}{2}\right)
    \mathbf{diag}(\beta)^{-1}$$
 \item moments of $V$:
 $$\mathbf{n}_{j,\cdot}^t :=\frac{2\alpha}{n} \mathcal{W}_{j}
\sum_{k:j_k=j} Y_{i_k,j_k} \mathbf{m}_{i_k,\cdot}^t $$
$$
 \mathcal{W}_j^{-1} := \frac{2\alpha}{n} \sum_{k:j_k=j}
  \left[\mathcal{V}_{i_k} + \mathbf{m}_{i_k,\cdot} \mathbf{m}_{i_k,\cdot}^t \right]
    + \left(a+\frac{m_1+m_2}{2}\right)
    \mathbf{diag}(\beta)^{-1}$$
 \item parameter of $\gamma$:
 $$\beta_k:= \frac{1}{2}\left[
 \sum_{i=1}^{m_1}\left(\mathbf{m}_{i,k}^2 + (\mathcal{V}_i)_{k,k} \right)
 +
 \sum_{j=1}^{m_2}\left(\mathbf{n}_{j,k}^2 + (\mathcal{V}_j)_{k,k} \right)
 \right]$$
\end{enumerate}
(where
$(\mathcal{V}_i)_{k,k}$ denotes the $(k,k)$-th entry of the matrix $\mathcal{V}_i$
and
$(\mathcal{W}_j)_{k,k}$ denotes the $(k,k)$-th entry of the matrix $\mathcal{W}_j$).

\subsection{Concentration of the posteriors}

For $r\geq 1$ and $B>0$ we define $\mathcal{M}(r,B)$ as the set of pairs of matrices $(\bar{U},\bar{V})$ with dimensions $m\times K$ and $p\times K$ respectively, satisfying the following constraints: $\bar{U}_{i,\ell}=0$ for $i>r$, $\|\bar{U}\|_{\infty} := \max_{i,\ell} |\bar{U}_{i,\ell}| \leq B$ and similarly $\bar{V}_{j,\ell}=0$ for $j>r$ and $\|\bar{V}\|_{\infty} \leq B$.
\begin{thm} Fix $a$ as any constant. There is a small enough $b>0$ such that
 \label{thm-matrix}
\begin{multline*}
\mathbb{E} \left[ \int d_{\alpha,\sigma}(M,M_0) \tilde{\pi}_{n,\alpha}({\rm d}M|X_1^n) \right]
\\
\leq \inf_{1\leq r \leq K} \inf_{(\bar{U},\bar{V})\in\mathcal{M}(r,B)}  \Biggl\{ \frac{\alpha}{1-\alpha} \frac{\left[\|M_0-\bar{U}\bar{V}^t\|_F + \frac{\sqrt{B}}{n}\right]^2}{2\sigma^2 mp}
 \\
 + \frac{2(1+\alpha)(1+2a) r(m+p) \left[ \log(nmp)
  + \mathcal{C}(a)
  \right]}{n(1-\alpha)} \Biggr\}
\end{multline*}
where the constant $\mathcal{C}(a)=\log(8\sqrt{\pi}\Gamma(a)2^{10a+1})+3$. In particular, the result holds for the choice $b=B^2/\{512 (nmp)^4 [(m\vee p) K]^2\}$.
\end{thm}
In practice, it is important that $b$ is small to ensure a good approximation of low-rank matrices~\cite{alquier2014bayesian}. We don't claim that $b= B^2 / \{512 (nmp)^4 [(m\vee p) K]^2\}$ is the optimal value,~\cite{alquier2014bayesian} recommends cross-validation to tune $b$.

Note as a special case that when $M=\bar{U}\bar{V}^t$ for $(\bar{U},\bar{V})\in \mathcal{M}(r,B)$ then we have exactly
\begin{multline}
\label{pre-coro-matrix}
\mathbb{E} \left[ \int d_{\alpha,\sigma}(M,M_0) \tilde{\pi}_{n,\alpha}({\rm d}M|X_1^n) \right]
 \\
 \leq
  \frac{2(1+\alpha)(1+2a) r(m+p) \left[ \log(nmp)
  + \mathcal{C}(a) + \frac{\alpha B}{2\sigma^2 mp}
  \right]}{n(1-\alpha)}.
\end{multline}
This result is the first consistency result for the VB approximation with Gaussian priors, that is used in practice. Still, it is stated for a ``weak'' distance criterion $d_{\alpha,\sigma}(M,M_0)$. Under additional assumptions, it is actually possible to relate this criterion to the standard Frobenius norm. Assume that there is a known $C$ such that $\max_{i,j}|(M_0)_{i,j}| \leq C$. This assumption is satisfied in many applications like collaborative filtering: in the Netflix data the entries are between $1$ and $5$. Then it is natural to project any estimator to the set of matrices with bounded entries. Precisely, define for any $M$ the matrix ${\rm clip}_C(M)$ its $(i,j)$-th entry: $\min(\max(M_{i,j},-C),C)$. A simple study of $d_{\alpha,\sigma}$, detailed in the proofs section, leads to the following result.
\begin{cor}
\label{cor-matrix}
 Under the assumptions of Theorem~\ref{thm-matrix}, and when in addition $\max_{i,j}|(M_0)_{i,j}| \leq C$, then
 \begin{multline*}
\mathbb{E} \left[ \int \|{\rm clip}_C(M) - M_0 \|_F^2 \tilde{\pi}_{n,\alpha}({\rm d}M|X_1^n) \right]
 \\
 \leq
  \frac{8C^2 (1+\alpha)(1+2a) r(m+p) \left[ \log(nmp)
  + \mathcal{C}(a) + \frac{\alpha B}{2\sigma^2 mp}
  \right]}{n[1-\exp(2 C^2 \alpha(\alpha-1) / \sigma^2)]}.
 \end{multline*}
\end{cor}
Note that once the Gaussian approximation of the posterior is known, it is easy to sample from it and to clip the samples to approximate the posterior mean of ${\rm clip}(M)$. So under the boundedness assumption we have a bound based on the Frobenius norm for an effective procedure based on VB. It is known that for the squared Frobenius norm, the rate $r(m+p)/n$ is minimax optimal -- maybe up to log terms~\cite{koltchinskii2011nuclear}.

Still assuming that $M=\bar{U}\bar{V}^t$ for $(\bar{U},\bar{V})\in \mathcal{M}(r,B)$ it is also possible to state a proper concentration result as an application of Corollary~\ref{cor-concentration}. We omit the proof as it is exactly similar to the one of Theorem~\ref{thm-matrix}.
\begin{thm}
Assume $M=\bar{U}\bar{V}^t$ for $(\bar{U},\bar{V})\in \mathcal{M}(r,B)$ and take $b$ as in Theorem~\ref{thm-matrix}. Then
  $$
  \mathbb{P}\left[
\int d_{\alpha,\sigma}(M,M_0) \tilde{\pi}_{n,\alpha}({\rm d}M|X_1^n)  \leq
 \frac{2(\alpha+1)}{1-\alpha} \varepsilon_n
 \right] \geq 1-\frac{2}{n\varepsilon_n}
 $$
 where for some explicit constant $\mathcal{D}(a,\sigma^2,B)$,
 $$\varepsilon_n = \frac{\mathcal{D}(a,\sigma^2,B)  r(m+p)\log(nmp)}{n} .$$
\end{thm}

\section{Nonparametric regression estimation}
\label{sec:density}

In this section, we provide a nonparametric example. Thus, the parameter will actually be a function $f$. We assume that $X_1=(W_1,Y_1),\dots,X_n=(W_n,Y_n)$ are i.i.d from a distribution $P_{f_0}$, and the model $(P_f)$ is given by: $W_i\sim \mathcal{U}([-1,1])$ and
$$ Y_i = f(W_i)+\xi_i  $$
where $\xi_i \sim \mathcal{N}(0,1)$. We will provide a prior and a mean-field approximation of the posterior. We will show that we estimate the functions $f$ belonging to a Sobolev ellipsoid $\mathcal{W}\left(r,C^2\right)$ at the minimax rate of convergence, up to $\log$ terms (the definitions of the ellipsoids will be reminded below). The reader might think that this example is not the most striking application of VB. On the other hand, it is an illustration of the generality of our method. We will estimate $f$ using projections on the Fourier basis and the choice of the number of coefficients will be done by model selection. It appears that in this case, model selection can be seen as a variational approximation where the constraint on the posterior is to give all its mass to only one model. This leads to adaptation of the estimator, in the sense that it is not required to know $r$ nor $C$ to compute the estimator.

\subsection{Construction of the prior}

First, we recall the definition of the trigonometric basis $(\varphi_k)_{k=1}^{\infty}$:
$$\varphi_1(t)=1,\, \varphi_{2k}(t)=\cos(\pi k t),\,
\varphi_{2k+1}(t)= \sin (\pi k t),\quad k=1, 2, \hdots $$
We now define a prior distribution $\pi$ by describing how to draw from $\pi$: we first draw $K$ from a geometric distribution, $\pi(K=k)=2^{-k}$. We then draw $\beta_1,\dots,\beta_K$ i.i.d from a $\mathcal{N}(0,1)$ distribution. We finally put
$$ f(x) = \sum_{k=1}^K \beta_k \varphi_k(x). $$
Note that when $f_0 (\cdot) = \sum_{k=1}^\infty \beta_k^0 \varphi_k(\cdot) $ and all the $\beta_k^0$'s are non-zero, such a function is never ``produced'' by the prior. Still, we will see that the prior gives enough mass to functions in the neighborhood of $f_0$, ensuring consistency.

\subsection{Construction of the variational approximation}

Note that the support of $\pi(\cdot|K)$ has dimension $K$, but the support of $\pi$ is infinite-dimensional. Thus, we can expect the support of the tempered posterior $\pi_{n,\alpha}$ to be also infinite-dimensional, and $\pi_{n,\alpha}$ to be intractable. We define a variational approximation that will fix these problems.

First, for $K\geq 1$ define $\mathcal{F}_K$ as the set of probability measures $\rho_{\mathbf{m},s^2}$ where $\mathbf{m}=(m_1,\dots,m_K)$ on functions $f(\cdot) = \sum_{k=1}^K \beta_k \varphi_k(\cdot)$
such that under $\rho_{\mathbf{m},s^2}$, the $\beta_k$'s are independent and $\beta_k \sim\mathcal{N}(m_k,s^2)$. We put $\mathcal{F} = \bigcup_{k=1}^{\infty} \mathcal{F}_k$. Note that the choice of a constant variance $s^2$ was motivated by the fact that the estimator of $\beta_k$ studied for example in~\cite{tsybakov2009book}, $\hat{\beta}_k = (1/n) \sum_{i=1}^n Y_i \varphi_k(X_i)$, satisfies $\hat{\beta}_k\sim\mathcal{N}(\beta_k,\sigma^2/n)$. Then
\begin{align*}
 \tilde{\pi}_{n,\alpha}
 & = \argmin_{\rho\in \mathcal{F}} \left\{ \alpha \int r_n(f,f_0) \rho({\rm d}f) + \mathcal{K}(\rho,\pi) \right\}
 \\
 & = \argmin_{K\geq 1} \argmin_{\mathbf{m},s^2} \Biggl\{ \alpha \int \frac{1}{2} \sum_{i=1}^n \left(Y_i - \sum_{k=1}^K \beta_k \varphi_k(W_i) \right)^2 \Phi({\rm d}\beta;\mathbf{m},s^2 I)
 \\
 & \quad \quad + \sum_{k=1}^K \frac{1}{2} \left[\log\left(\frac{1}{s^2}\right)+s^2+m_k^2-1\right]  + K\log(2) \Biggr\}
 \\
 & = \argmin_{K\geq 1} \argmin_{\mathbf{m},s^2} \Biggl\{   \frac{\alpha}{2} \sum_{i=1}^n \left(Y_i - \sum_{k=1}^K m_k \varphi_k(W_i) \right)^2 + \frac{s^2 \alpha}{2}\sum_{i=1}^n \sum_{k=1}^K \varphi^2_k(W_i)
 \\
 & \quad \quad + \sum_{k=1}^K \frac{1}{2} \left[\log\left(\frac{1}{s^2}\right)+s^2+m_k^2-1\right]  + K\log(2) \Biggr\}
\end{align*}
(that is, the approximated posterior mean is simply a ridge regression estimator).

\subsection{Nonparametric rates of convergence}

We remind the definition of the Sobolev ellipsoid given (see {\it e.g.} Chapter 1 in~\cite{tsybakov2009book}) for $C>0$ and $r\geq 2$:
$$\mathcal{W}\left(r,C^2\right)= \left\{ f\in L_{2}([-1,1]):
f=\sum_{k=1}^{\infty}\beta_k \varphi_k \mbox{ and } 
     \sum_{k=1}^{\infty} k^{2r}
\beta_{k}^{2}  \leq C^2 \right\}.
$$
\begin{thm}
\label{thm-concentration-vb:expect:regression}
Fix $\alpha\in(0,1)$. Assume that there is an $r\in[2,\infty[$ and a $C>0$ such that $f_0\in\mathcal{W}(r,C^2)$. Then
\begin{equation*}
\mathbb{E}  \left[ \int D_{\alpha}(P_{f},P_{f_0}) \tilde{\pi}_{n,\alpha}({\rm d}\theta|X_1^n) \right] =\mathcal{O} \left( \left(\frac{\log(n)}{n}\right)^{\frac{2r}{2r+1}}  \right).
\end{equation*}
\end{thm}
The proof is in Section~\ref{section:proofs}. Note that on the contrary to previous sections, we only provide an asymptotic statement here. However, from the proof of Theorem~\ref{thm-concentration-vb:expect:regression}, it is clear that it is possible to provide a non-asymptotic statement as well (with cumbersome constants).

Here again, note that the distance criterion used in the left-hand side is not standard. We actually have:
$$ D_{\alpha}(P_{f},P_{f_0}) = \frac{1}{\alpha-1} \log \left[ \frac{1}{2} \int_{-1}^1 \exp\left(\frac{\alpha(\alpha-1)(f(x)-f_0(x))^2}{2}\right) {\rm d} x \right] . $$
However, when $f_0 \in\mathcal{W}_{r,C^2}$, $f$ is bounded by a constant that depends on $r$ and $C^2$. If we moreover assume that $f_0$ is bounded by a known constant $c_0$, we can as in Section~\ref{sec:matrix-compl} define a clip operator: ${\rm clip}_{c_0}(f)(x) = \min(\max(-c_0,f(x)),c_0)$ and obtain:
\begin{equation*}
\mathbb{E}  \left[ \int  \left\|{\rm clip}_{c_0}(f) - f_0 \right\|_2^2 \tilde{\pi}_{n,\alpha}({\rm d}\theta|X_1^n) \right] =\mathcal{O} \left( \left(\frac{\log(n)}{n}\right)^{\frac{2r}{2r+1}}  \right).
\end{equation*}
The rate $1/n^{2r/(2r+1)}$ is known to be minimax optimal on $\mathcal{W}(r,C^2)$ for the squared  $\|\cdot\|_2$-norm~\cite{tsybakov2009book}. The additional $\log$ term is sometimes referred to as ``the price to pay for adaptation''. In the case of the $\|\cdot\|_2$-norm this is misleading as it is actually possible to build an adaptive estimator that reaches the minimax rate without the additional $\log$, but up to our knowledge this is not possible with a fully Bayesian estimator.

\section{Conclusion}

Based on PAC-Bayesian inequalities, we introduced a generic method to study the concentration of variational Bayesian approximations. This is a very general approach that can be applied to many models. We studied applications to logistic regression, matrix completion and density estimation. Still, some questions remain open. From a theoretical perspective, the oracle inequality in Theorem~\ref{thm-concentration-vb:expect:misspecified} compares a R\'enyi divergence to a Kullback-Leibler divergence. It would be very interesting to obtain a result with the Kullback divergence in the left-hand side. This is probably more difficult, if possible at all. We believe that tools from~\cite{grunwald2016fast} could be of some help, but some work is needed to make explicit the assumptions of this paper in our context.

Also, since the first version of this work was submitted, extensions were proven by other authors:~\cite{yangnew} extended our results to models with hidden variables, such as mixture models, and~\cite{gaonew} proved results in the case $\alpha=1$ and study many nonparametric examples. Note that while $\alpha=1$ remains the most popular choice in practice, these results require much stronger assumptions and cannot in general be extended to the misspecified case~\cite{grunwaldNON}.

An important open issue is the choice of the parameter $\alpha$. It is clear that our results are not helpful to solve this issue. Some previous work proposes to use cross-validation~\cite{alquier2016properties}, but this is computationaly expensive. Moreover, no theoretical guarantees are known in this case. In the misspecified case,~\cite{grunwald2012safe} proposed an online adaptive tuning of this parameter. However, it is not clear if this method could work in our context. This should be the object of a future work.

Finally, it would be nice to get rid of the extra log in the rates. Catoni's localization technique~\cite{MR2483528} is a nice tool to remove extra log factors in PAC-Bayesian bounds, but its adaptation to our setting is not direct. It could be the object of future works.

\section*{Acknowledgements}

We would like to thank Badr-Eddine Ch\'erief-Abdellatif, as well as the Associate Editor and the anonymous Referees, for their helpful comments and suggestions on the paper.

\section{Proofs}
\label{section:proofs}
\subsection{Proof of Theorem~\ref{thm-bha}}

We adapt the proof given in~\cite{bhattacharya2016bayesian}. Fix $\alpha\in(0,1)$, and $\theta\in\Theta$. It's immediate to check that
\begin{equation*}
 \mathbb{E}\left[ \exp\left(-\alpha r_n(\theta,\theta_0)\right) \right] = \exp\left[-(1-\alpha) D_\alpha(P_{\theta}^{\otimes n},P_{\theta_0}^{\otimes n}) \right].
\end{equation*}
Note that it might be that $D_\alpha(P_{\theta}^{\otimes n},P_{\theta_0}^{\otimes n})=n D_\alpha(P_{\theta},P_{\theta_0}) = +\infty$. Rewrite $\pi({\rm d}\theta) = \pi({\rm d}\theta)\mathbf{1}_{\{D_\alpha(P_{\theta},P_{\theta_0}) < +\infty\}} + \pi({\rm d}\theta)\mathbf{1}_{\{D_\alpha(P_{\theta},P_{\theta_0}) = +\infty\}} =\pi_1({\rm d}\theta) + \pi_2({\rm d}\theta)$. First, when $\pi=\pi_2$, $\pi$-almost surely, $P_\theta$ is singular to $P_{\theta_0}$. But then, $\pi$-almost surely, $D_\alpha(P_\theta,P_{\theta_0})=+\infty$ and $r_n(\theta,\theta_0)=+\infty$. This also holds $\rho$-almost surely for any $\rho \ll \pi$ and thus the statement of the theorem is trivial: $\mathbb{P}(+\infty\leq +\infty)\geq 1-\varepsilon$.

Assume now that $\pi\neq \pi_2$. This allows to define the renormalization $\tilde{\pi}(\cdot)= \pi_1(\cdot) / \pi(D_\alpha(P_{\theta},P_{\theta_0}) < +\infty) $, that is a probability measure. On the support of $\tilde{\pi}(\cdot)$,
\begin{equation*}
 \mathbb{E}\left[ \exp\left(-\alpha r_n(\theta,\theta_0) + (1-\alpha)n D_{\alpha}(P_{\theta},P_{\theta_0}) \right) \right] = 1.
\end{equation*}
Integrate with respect to $\tilde{\pi}$,
\begin{equation*}
 \int \mathbb{E}\left[ \exp\left(-\alpha r_n(\theta,\theta_0) + (1-\alpha)n D_{\alpha}(P_{\theta},P_{\theta_0})  \right) \right] \tilde{\pi}({\rm d}\theta) = 1
\end{equation*}
and using Fubini's theorem,
\begin{equation}
\label{proofcaca}
 \mathbb{E}\left[ \int \exp\left(-\alpha r_n(\theta,\theta_0) + (1-\alpha)n D_{\alpha}(P_{\theta},P_{\theta_0})  \right) \tilde{\pi}({\rm d}\theta) \right]  = 1.
\end{equation}
The key argument here, introduced by~\cite{catoni2004statistical}, is to use Lemma~\ref{thm-dv}. Note that almost surely with respect to the sample, we know that
$$h(\theta) := - \alpha r_n(\theta,\theta_0) + (1-\alpha)n D_{\alpha}(P_{\theta},P_{\theta_0})$$
satisfies $ \int \exp(h) {\rm d}\tilde{\pi} < \infty$, otherwise, the expectation in~\eqref{proofcaca} would be infinite. So, the conditions of Lemma~\ref{thm-dv} are satisfied almost surely with respect to the sample, and we obtain:
\begin{multline*}
 \mathbb{E}\Biggl\{  \exp\Biggl[ \sup_{\rho\in\mathcal{M}_{1}^+(\Theta)} \biggl( \int\left(-\alpha r_n(\theta,\theta_0) + (1-\alpha)n D_{\alpha}(P_{\theta},P_{\theta_0})\right)\rho({\rm d}\theta) 
 \\
 -\mathcal{K}(\rho,\tilde{\pi})\biggr) \Biggr] \Biggr\}  = 1.
\end{multline*}
Multiply both sides by $\varepsilon$ to get
\begin{multline*}
 \mathbb{E}\Biggl\{  \exp\Biggl[ \sup_{\rho\in\mathcal{M}_{1}^+(\Theta)} \biggl( \int\left(-\alpha r_n(\theta,\theta_0) + (1-\alpha)n D_{\alpha}(P_{\theta},P_{\theta_0})\right)\rho({\rm d}\theta) 
 \\
 -\mathcal{K}(\rho,\tilde{\pi})\biggr) - \log\left(\frac{1}{\varepsilon} \right)\Biggr] \Biggr\}  = \varepsilon.
\end{multline*}
Using Markov's inequality,
\begin{multline*}
 \mathbb{P}\Biggl[ \sup_{\rho\in\mathcal{M}_{1}^+(\Theta)}\biggl(   \int\left(-\alpha r_n(\theta,\theta_0) + (1-\alpha)n D_{\alpha}(P_{\theta},P_{\theta_0})\right)\rho({\rm d}\theta) 
 \\
 -\mathcal{K}(\rho,\tilde{\pi})\biggr)- \log\left(\frac{1}{\varepsilon} \right)  \geq 0 \Biggr] \leq \varepsilon.
\end{multline*}
Taking the complementary event,
\begin{multline*}
 \mathbb{P}\Biggl( \forall \rho\in\mathcal{M}_{1}^+(\Theta),\quad \int\left(-\alpha r_n(\theta,\theta_0) + (1-\alpha)n D_{\alpha}(P_{\theta},P_{\theta_0})\right)\rho({\rm d}\theta) 
 \\
 -\mathcal{K}(\rho,\tilde{\pi})- \log\left(\frac{1}{\varepsilon} \right)  \leq 0 \Biggr) \geq 1- \varepsilon.
\end{multline*}
Now, for a given $\rho$, it might be that
$\int n D_{\alpha}(P_{\theta},P_{\theta_0})\rho({\rm d}\theta) =\infty $
but then, the previous equation implies that
$ \int r_n(\theta,\theta_0) \rho({\rm d}\theta)  + \mathcal{K}(\rho,\tilde{\pi}) =\infty $ and so the statement of the theorem is trivially satisfied as $\infty \leq \infty$. On the other hand, assuming that $\int n D_{\alpha}(P_{\theta},P_{\theta_0})\rho({\rm d}\theta) < \infty $
we rearrange terms to get
 \begin{multline*}
 \mathbb{P}\Biggl( \forall \rho\in \mathcal{M}_{1}^+(\Theta),
 \int  D_{\alpha}(P_{\theta},P_{\theta_0}) \rho({\rm d}\theta) \Biggr.
 \\
 \Biggl.\leq \frac{\alpha}{1-\alpha} \int \frac{r_n(\theta,\theta_0) }{n} \rho({\rm d}\theta)
 + \frac{\mathcal{K}(\rho,\tilde{\pi}) + \log\left(\frac{1}{\varepsilon}\right)}{n(1-\alpha)}
 \Biggr) \geq 1-\varepsilon.
 \end{multline*}
Now, we decompose $\rho = \rho_1 + \rho_2$ as we decomposed $\pi$. First, when $\rho \neq \rho_1$, we have: $\rho(D_{\alpha}(P_{\theta},P_{\theta_0}) =+\infty) >0 $. But then this means that $\rho( r_n(\theta,\theta_0) = +\infty)>0$, and once again, the statement of the theorem is trivial: $\mathbb{P}(+\infty\leq +\infty)\geq 1-\varepsilon$. So we can assume that $\rho = \rho_1$. But then $\mathcal{K}(\rho,\tilde{\pi}) = \mathcal{K}(\rho,\pi) + \log \pi\left(D_{\alpha}(P_{\theta},P_{\theta_0}) <+\infty \right) \leq \mathcal{K}(\rho,\pi)  $, thus the statement of the theorem also holds. This ends the proof.

\subsection{Proof of Theorem~\ref{thm-concentration-vb}}

Fix $\eta\in(0,1)$ and define
\begin{multline*}
\rho^* = \underset{\rho \in \mathcal{F}}{\arg\min} \, \Biggl\{ \frac{\alpha}{1-\alpha} \int \frac{ \mathbb{E}[r_n(\theta,\theta_0)]}{n}\rho({\rm d}\theta)
\\
+ \frac{\alpha}{n(1-\alpha)}\sqrt{\frac{{\rm Var}[\int r_n(\theta,\theta_0)\rho({\rm d}\theta)] }{\eta}} + \frac{\mathcal{K}(\rho,\pi)}{n(1-\alpha)} \Biggr\}.
\end{multline*}
Chebyshev's inequality leads to
\begin{multline*}
\mathbb{P}\Biggl\{
\frac{\alpha}{1-\alpha} \int \frac{r_n(\theta,\theta_0)}{n} \rho^*({\rm d}\theta) \geq 
\frac{\alpha}{1-\alpha} \int \frac{ \mathbb{E}[r_n(\theta,\theta_0)]}{n}\rho^*({\rm d}\theta) \\
+ \frac{\alpha}{n(1-\alpha)}\sqrt{\frac{{\rm Var}[\int r_n(\theta,\theta_0)\rho^*({\rm d}\theta)] }{\eta}} + \frac{\mathcal{K}(\rho^*,\pi)}{n(1-\alpha)} 
\Biggr\} \leq \eta
\end{multline*}
and so
\begin{multline}
\label{step:cheby}
\mathbb{P}\Biggl\{
\frac{\alpha}{1-\alpha} \int \frac{r_n(\theta,\theta_0)}{n} \rho^*({\rm d}\theta) \geq 
\frac{\alpha}{1-\alpha} \int \mathcal{K}(P_{\theta_0},P_{\theta}) \rho^*({\rm d}\theta) \\
+ \frac{\alpha}{1-\alpha}\sqrt{ \frac{\int{\rm Var}\left[\log\frac{p_{\theta}(X_i)}{p_{\theta_0}(X_i)}\right] \rho^*({\rm d}\theta) }{n \eta}} + \frac{\mathcal{K}(\rho^*,\pi)}{n(1-\alpha)}
\Biggr\} \leq \eta.
\end{multline}
Now apply take the union bound of this inequality and of the inequality in Corollary \ref{thm-cor-bha}. We obtain, for any $\alpha\in(0,1)$, for any $\varepsilon\in(0,1)$, with probability at least $1-\varepsilon-\eta$,
 \begin{align*}
 & \int D_{\alpha}(P_{\theta},P_{\theta_0}) \pi_{n,\alpha}({\rm d}\theta|X_1^n) \\
 & \quad \leq  \inf_{\rho\in\mathcal{F}} \left\{ \frac{\alpha \int r_n(\theta,\theta_0) \rho({\rm d}\theta) + \mathcal{K}(\rho,\pi) + \log\left(\frac{1}{\varepsilon}\right)}{1-\alpha} \right\} \text{ by Cor. \ref{thm-cor-bha}} \\
 & \quad \leq    \frac{\alpha \int r_n(\theta,\theta_0) \rho^*({\rm d}\theta) + \mathcal{K}(\rho^*,\pi) + \log\left(\frac{1}{\varepsilon}\right)}{1-\alpha} \\
 & \quad \leq  \frac{\alpha \int \left[ \mathcal{K}(P_{\theta_0},P_{\theta})  + \sqrt{ \frac{1}{n\eta}{\rm Var}[\int r_n(\theta,\theta_0)\rho^*({\rm d}\theta)]} \right]  }{1-\alpha} \\
 & \quad \quad \quad + \frac{\mathcal{K}(\rho^*,\pi) + \log\left(\frac{1}{\varepsilon}\right)}{n(1-\alpha)} \text{ by~\eqref{step:cheby}} \\
 & \quad = \inf_{\rho\in\mathcal{F}} \Biggl\{ \frac{\alpha \int \left[ \mathcal{K}(P_{\theta_0},P_{\theta})  + \sqrt{ \frac{1}{n\eta}{\rm Var}[\int r_n(\theta,\theta_0)]} \right] \rho({\rm d}\theta) }{1-\alpha}
 \\
 & \quad \quad \quad + \frac{\mathcal{K}(\rho,\pi) + \log\left(\frac{1}{\varepsilon}\right)}{n(1-\alpha)}
 \Biggr\} \text{ by definition of } \rho^*
 \\
  & \quad \leq \inf_{\rho\in\mathcal{F}} \Biggl\{ \frac{\alpha \int \left[ \mathcal{K}(P_{\theta_0},P_{\theta})  + \sqrt{ \frac{1}{n\eta}\mathbb{E} \left[ \log^2 \left( \frac{p_{\theta}(X_i)}{p_{\theta_0}(X_i)} \right) \right]} \right] \rho({\rm d}\theta) }{1-\alpha}
 \\
 & \quad \quad \quad + \frac{\mathcal{K}(\rho,\pi) + \log\left(\frac{1}{\varepsilon}\right)}{n(1-\alpha)}
 \Biggr\}
 \\
 & \quad \leq  \frac{\alpha \int \left[ \mathcal{K}(P_{\theta_0},P_{\theta})  + \sqrt{ \frac{1}{n\eta}\mathbb{E} \left[ \log^2 \left( \frac{p_{\theta}(X_i)}{p_{\theta_0}(X_i)} \right) \right]} \right] \rho_n({\rm d}\theta) }{1-\alpha}
 \\
 & \quad \quad \quad + \frac{\mathcal{K}(\rho_n,\pi) + \log\left(\frac{1}{\varepsilon}\right)}{n(1-\alpha)}
 \\
 & \quad \leq \frac{\alpha\left(\varepsilon_n + \sqrt{\frac{\varepsilon_n}{n\eta}}\right)}{1-\alpha}
             + \frac{n \varepsilon_n + \log\left(\frac{1}{\varepsilon}\right)}{n(1-\alpha)}
 \end{align*}
 where in the last step we use the assumptions on $\rho_n$.

\subsection{Proof of Theorem~\ref{thm-concentration-vb:expect}}

The beginning is as for Theorem~\ref{thm-bha}. Fix $\alpha\in(0,1)$, then
\begin{equation*}
 \mathbb{E}\left[ \exp\left(-\alpha r_n(\theta,\theta_0) - (1-\alpha)D_{\alpha}(P_{\theta}^{\otimes n},P_{\theta_0}^{\otimes n}) \right) \right] = 1.
\end{equation*}
Integrate with respect to $\pi$,
\begin{equation*}
 \int \mathbb{E}\left[ \exp\left(-\alpha r_n(\theta,\theta_0) - (1-\alpha)D_{\alpha}(P_{\theta}^{\otimes n},P_{\theta_0}^{\otimes n})\right) \right] \pi({\rm d}\theta) = 1
\end{equation*}
and using Fubini's theorem and Lemma~\ref{thm-dv}
\begin{multline*}
 \mathbb{E}\Biggl\{  \exp\Biggl[ \sup_{\rho\in\mathcal{M}_{1}^+(\Theta)} \biggl( \int\left(-\alpha r_n(\theta,\theta_0) - (1-\alpha)D_{\alpha}(P_{\theta}^{\otimes n},P_{\theta_0}^{\otimes n})\right)\rho({\rm d}\theta) 
 \\
 -\mathcal{K}(\rho,\pi)\biggr)\ \Biggr] \Biggr\}  = 1.
\end{multline*}
This is where things change: we now use Jensen's inequality to obtain
\begin{multline*}
\mathbb{E}\Biggl[ \sup_{\rho\in\mathcal{M}_{1}^+(\Theta)} \biggl( \int\left(-\alpha r_n(\theta,\theta_0) - (1-\alpha)D_{\alpha}(P_{\theta}^{\otimes n},P_{\theta_0}^{\otimes n})\right)\rho({\rm d}\theta) 
 \\
 -\mathcal{K}(\rho,\pi)\biggr) \Biggr]  = 0
\end{multline*}
and so as a special case
\begin{multline*}
\mathbb{E}\Biggl[  \int\left(-\alpha r_n(\theta,\theta_0) - (1-\alpha)D_{\alpha}(P_{\theta}^{\otimes n},P_{\theta_0}^{\otimes n})\right)\tilde{\pi}_{n,\alpha}({\rm d}\theta|X_1^n)
 \\
 -\mathcal{K}(\tilde{\pi}_{n,\alpha}(\cdot|X_1^n),\pi) \Biggr]  = 0.
\end{multline*}
Rearranging terms,
\begin{align*}
\mathbb{E} & \left[  \int D_{\alpha}(P_{\theta}^{\otimes n},P_{\theta_0}^{\otimes n})\tilde{\pi}_{n,\alpha}({\rm d}\theta|X_1^n) \right]
\\
& \leq  \mathbb{E} \left[ \frac{ \alpha}{1-\alpha} \int r_n(\theta,\theta_0) \tilde{\pi}_{n,\alpha}({\rm d}\theta|X_1^n)  + \frac{\mathcal{K}(\tilde{\pi}_{n,\alpha}(\cdot|X_1^n),\pi) }{1-\alpha}
\right]
\\
& = \mathbb{E} \left\{ \inf_{\rho\in\mathcal{F}}  \left[\frac{ \alpha}{1-\alpha} \int r_n(\theta,\theta_0) \rho({\rm d}\theta)  + \frac{\mathcal{K}(\rho,\pi) }{1-\alpha}
\right]\right\} \text{ by dfn.}
\\
& \leq \inf_{\rho\in\mathcal{F}} \left\{  \mathbb{E}  \left[\frac{ \alpha}{1-\alpha} \int r_n(\theta,\theta_0) \rho({\rm d}\theta)   + \frac{\mathcal{K}(\rho,\pi) }{1-\alpha}
\right]\right\}
\\
& = \inf_{\rho\in\mathcal{F}} \left\{ \frac{n \alpha}{1-\alpha} \int \mathcal{K}(P_{\theta_0},P_{\theta}) \rho({\rm d}\theta)  + \frac{\mathcal{K}(\rho,\pi) }{1-\alpha}\right\}
\end{align*}
and so
\begin{align*}
\mathbb{E}  \left[  \int D_{\alpha}(P_{\theta},P_{\theta_0})\tilde{\pi}_{n,\alpha}({\rm d}\theta|X_1^n) \right]
& =
\mathbb{E} \left[  \int \frac{D_{\alpha}(P_{\theta}^{\otimes n},P_{\theta_0}^{\otimes n})}{n}\tilde{\pi}_{n,\alpha}({\rm d}\theta|X_1^n) \right]
\\
& \leq \inf_{\rho\in\mathcal{F}} \left\{ \frac{ \alpha}{1-\alpha} \int \mathcal{K}(P_{\theta_0},P_{\theta}) \rho({\rm d}\theta)  + \frac{\mathcal{K}(\rho,\pi) }{n(1-\alpha)}\right\}.
\end{align*}

\subsection{Proof of Theorem \ref{thm-lipschitz-dens}}

We start by defining a sequence $\rho_n(d\theta):=\Phi(d\theta;\theta_0,\sigma_n^2I)\in\mathcal{F}_\Phi^{id}$ indexed by a positive scalar $\sigma^2_n$ to be later defined. As before by proving the result on the smallest family of distribution, it will remain true on larger ones using the fact that $\min_{\mathcal{F}^{id}}\leq \min_{\mathcal{F}^{diag}}\leq \min_{\mathcal{F}^{full}}$.
Under Assumption \ref{lipschitz-den} we can check the hypotheses on the KL between the likelihood terms as required in Theorem \ref{thm-concentration-vb}.
We have 
\[
\mathcal{K}(P_{\theta_0},P_\theta)= \mathbb{E}\left[\log p_{\theta_0}(X)-\log p_\theta(X)\right]\leq \mathbb{E}\left[M(X)\right]\Vert\theta-\theta_0\Vert_2
\]
and 
\[
\mathbb{E}\left[\log^2\frac{p_{\theta_0}}{p_\theta}(X)\right]= \mathbb{E}\left[\left(\log p_{\theta_0}(X)-\log p_\theta(X)\right)^2\right]\leq \mathbb{E}\left[M(X)^2\right]\Vert\theta-\theta_0\Vert_2
\]
When integrating with respect to $\rho_n$ we have 
$$
  \int\mathcal{K}(P_{\theta_0},P_{\theta})\rho_n(d\theta)\leq B_1\sigma_n\sqrt{d} \text{ and }
 \int\mathbb{E}\left[\log^2\frac{p_{\theta_0}}{p_\theta}(X)\right]\rho_n(d\theta)\leq B_2\sigma^2_nd.
 $$

To apply Theorem \ref{thm-concentration-vb} it remains to compute the KL between the approximation of the pseudo-posterior and the prior,
\begin{align*}
  \frac1n \mathcal{K}(\rho_n,\pi)&=\frac dn \left[
  \frac{1}{2}\log\left(\frac{\vartheta^2}{\sigma^2}\right)
        + \frac{\sigma^2}{\vartheta^2}
  \right]
  + \frac{\|\theta_0\|^2}{n\vartheta^2} - \frac{d}{2n}.
\end{align*}
To obtain an estimate of the rate $\varepsilon_n$ of Theorem \ref{thm-concentration-vb} we put together those bounds. Choosing $\sigma^{2}_n=\frac1{n\sqrt{d}}$ we can apply 
it with
\[
\varepsilon_n=\frac{B_1}n\vee \frac{B_2}{n^2} \vee \left\lbrace\frac dn \left[
  \frac12\log\left(\vartheta^2 n^2 \sqrt{d}\right)
        + \frac{1}{n\vartheta^2}
  \right]
  + \frac{\|\theta_0\|^2}{n\vartheta^2} - \frac{d}{2n}\right\rbrace.
\]

\subsection{Proof of Theorem \ref{thm-vb-convex}}

From the proof of Theorem  \ref{thm-concentration-vb:expect} we get 
\begin{align*}
\mathbb{E} & \left[  \int D_{\alpha}(P_{\theta}^{\otimes n},P_{\theta_0}^{\otimes n})\tilde{\pi}^k_{n,\alpha}({\rm d}\theta|X_1^n) \right]
\\
& \leq  \mathbb{E} \left[ \frac{ \alpha}{1-\alpha} \int r_n(\theta,\theta_0) \tilde{\pi}^k_{n,\alpha}({\rm d}\theta|X_1^n)  + \frac{\mathcal{K}(\tilde{\pi}^k_{n,\alpha}(\cdot|X_1^n),\pi) }{1-\alpha}
\right]
\\
& = \mathbb{E} \left\{ \inf_{\rho\in\mathcal{F}^\Phi_B}  \left[\frac{ \alpha}{1-\alpha} \int r_n(\theta,\theta_0) \rho({\rm d}\theta)  + \frac{\mathcal{K}(\rho,\pi) }{1-\alpha}
\right]\right\}\\ 
&\qquad+ \left\lbrace\mathbb{E} \left[ \frac{ \alpha}{1-\alpha} \int r_n(\theta,\theta_0) \tilde{\pi}^k_{n,\alpha}({\rm d}\theta|X_1^n)  + \frac{\mathcal{K}(\tilde{\pi}^k_{n,\alpha}(\cdot|X_1^n),\pi) }{1-\alpha}
\right]\right.\\
&\left.\qquad-\mathbb{E} \left\{ \inf_{\rho\in\mathcal{F}^\Phi_B}  \left[\frac{ \alpha}{1-\alpha} \int r_n(\theta,\theta_0) \rho({\rm d}\theta)  + \frac{\mathcal{K}(\rho,\pi) }{1-\alpha}
\right]\right\}\right\rbrace.
\end{align*}

By definition of $f$ we get
\begin{align*}
\mathbb{E} & \left[  \int D_{\alpha}(P_{\theta}^{\otimes n},P_{\theta_0}^{\otimes n})\tilde{\pi}^k_{n,\alpha}({\rm d}\theta|X_1^n) \right]\\
& = \mathbb{E} \left\{ \inf_{\rho\in\mathcal{F}^\Phi_B}  \left[\frac{ \alpha}{1-\alpha} \int r_n(\theta,\theta_0) \rho({\rm d}\theta)  + \frac{\mathcal{K}(\rho,\pi) }{1-\alpha}\right]\right\}
+\frac1{1-\alpha}\mathbb{E}\left\lbrace\mathbb{E}f(\bar{x}_k,\xi)-\inf_{u\in \mathbb{B}}\mathbb{E}f(u,\xi)\right\rbrace
\\
& \leq \inf_{\rho\in\mathcal{F}^\Phi_B} \left\{  \mathbb{E}  \left[\frac{ \alpha}{1-\alpha} \int r_n(\theta,\theta_0) \rho({\rm d}\theta)   + \frac{\mathcal{K}(\rho,\pi) }{1-\alpha}
\right]\right\}+\frac1{1-\alpha}\mathbb{E}\left\lbrace\mathbb{E}f(\bar{x}_k,\xi)-\inf_{u\in \mathbb{B}}\mathbb{E}f(u,\xi)\right\rbrace
\\
& = \inf_{\rho\in\mathcal{F}^\Phi_B} \left\{ \frac{n \alpha}{1-\alpha} \int \mathcal{K}(P_{\theta_0},P_{\theta}) \rho({\rm d}\theta)  + \frac{\mathcal{K}(\rho,\pi) }{1-\alpha}\right\}+\frac1{1-\alpha}\mathbb{E}\left\lbrace\mathbb{E}f(\bar{x}_k,\xi)-\inf_{u\in \mathbb{B}}\mathbb{E}f(u,\xi)\right\rbrace.
\end{align*}
Following the rest of the proof of \ref{thm-concentration-vb:expect} we get  
\begin{align*}
\mathbb{E} & \left[  \int D_{\alpha}(P_{\theta},P_{\theta_0})\tilde{\pi}^k_{n,\alpha}({\rm d}\theta|X_1^n) \right]\\
&\leq \inf_{\rho\in\mathcal{F}} \left\{ \frac{ \alpha}{1-\alpha} \int \mathcal{K}(P_{\theta_0},P_{\theta}) \rho({\rm d}\theta)  + \frac{\mathcal{K}(\rho,\pi) }{n(1-\alpha)}\right\}+\frac1{n(1-\alpha)}\mathbb{E}\left\lbrace\mathbb{E}f(\bar{x}_k,\xi)-\inf_{u\in \mathbb{B}}\mathbb{E}f(u,\xi)\right\rbrace.
\end{align*}

To bound the first term of the right hand-side we use Assumption \ref{lipschitz-den} and  
the proof of Theorem \ref{thm-lipschitz-dens}. In particular notice that $ \Phi(d\theta;\theta_0,\frac1{n\sqrt{d}}I_d)\in\mathcal{F}^\Phi_B$,
we get straight away
\begin{align*}
\mathbb{E} & \left[  \int D_{\alpha}(P_{\theta},P_{\theta_0})\tilde{\pi}^k_{n,\alpha}({\rm d}\theta|X_1^n) \right]
&\leq \frac{1+\alpha}{1-\alpha}\varepsilon_n+\frac1{n(1-\alpha)}\mathbb{E}\left\lbrace\mathbb{E}f(\bar{x}_k,\xi)-\inf_{u\in \mathbb{B}}\mathbb{E}f(u,\xi)\right\rbrace
\end{align*}
We now study the term inside the brackets on the right hand-side.

First notice that the sequence $(x_t)_{t\geq 0}$ in Algorithm \ref{algo-SVB} is equivalent to that of an online gradient descent on the sequence $\left\lbrace f(x,\xi_t)\right\rbrace_t$. Hence under Assumption \ref{algo-SVB} we can apply Corollary 2.7 of \cite{Shalev-Shwartz2012} with $\gamma_T=\frac{B}{L\sqrt{2T}}$ to get the following bound on the regret for any $u\in \mathbb{B}$ 
\[
\sum_{t=1}^T f(x_t,\xi_t)-\sum_{t=1}^T f(u,\xi_t)\leq \sqrt{2BLT}.
\] 
Divide by $T$, take expectation with respect to $(\xi_t)_t$  
\[
\frac1T \sum_{t=1}^T \mathbb{E}f(x_t,\xi_t)-\mathbb{E}f(u,\xi)\leq \sqrt{\frac{2BL}{T}}.
\] 
Notice that $x_t$ belongs to the $\sigma$-algebra generated by $(x_1,\cdots,x_{t-1})$. By a multiple use
of the tower property we get, 
\begin{align*}
\frac1T \sum_{t=1}^T \mathbb{E}f(x_t,\xi)-\mathbb{E}f(u,\xi)&\leq \sqrt{\frac{2BL}{T}}\\
\mathbb{E}f(\bar{x},\xi)-\mathbb{E}f(u,\xi)&\leq \sqrt{\frac{2BL}{T}},\text{ by Jensen and the convexity of } f.
\end{align*}
Putting everything together concludes the proof.
\subsection{Proof of Corollary \ref{cor:logit-cor}}
Direct calculation shows that the log-likelihood is $2\Vert X\Vert$-Lipschitz hence satisfying Assumption \ref{lipschitz-den}.
We conclude using Theorem \ref{thm-lipschitz-dens} and the assumption on the design matrix. 

\subsection{Proof of Corollary \ref{cor:optim-logit}}
Start by noticing that we can take $f$ as 
\[
f((m,C),\xi):=\alpha\log p_{m+C\xi}(\tilde{x})+\mathcal{K}(\rho,\pi),
\]
where $\rho(.)=\Phi(.;m,CC^t)$ the likelihood part is convex with Lipschitz gradient as a composition of a convex and gradient Lipschitz function with a 
affine map. The Lipschitz constant for this term is bounded by $\sum_{i=1}^n\Vert x_ix_i^t\Vert$. The KL part can be written as $\mathcal{K}(\rho,\pi)=\frac{\Vert m\Vert^2}{2\vartheta}+\left(\frac1{2\vartheta}\text{trace}(CC^t)-\log\vert C\vert\right)$ which is convex for positive semi-definite $C$. We need to check that the gradients of the objectives are also Lipschitz, the only problematic term is $\log \text{det}(C)$. Denote $(\lambda_i)$ the eigen values of $\Sigma=CC^t$ 
\begin{align*}
\Sigma\succeq\psi I_{d\times d}&\Rightarrow \forall i\in\lbrace1,\cdots,d\rbrace,\quad\frac1{\lambda_i}\leq \frac1\psi\\
&\Rightarrow \text{trace}(\Sigma^{-1})\leq \frac d{\psi}\\
&\Rightarrow \text{trace}^{\frac12}\left( C^{-1}C^{-T}\otimes C^{-1}C^{-T}\right)\leq \frac d{\psi}\\
&\Rightarrow\text{trace}^{\frac12}\left( (C^{-1}\otimes C^{-T})( C^{-1}\otimes C^{-T})\right)\leq \frac d{\psi}\\
&\Rightarrow\Vert\nabla^2_C\log\text{det} C\Vert_2\leq \frac d\psi.
\end{align*}
To apply Theorem \ref{thm-vb-convex} we also need to check that the new constraint contains the Gaussian distribution used in the proof. This is the case as long as $\psi\leq \sigma^2=\frac1{n\sqrt{d}}$.

\bigskip

\noindent \textbf{Supplemtary material}.
\begin{itemize}
\item The supplementary material contains the toy example mentioned in Remark~\ref{remark-laplace} above. 
\item The remaining proofs, that is, the proofs of Theorems~\ref{thm-matrix} and~\ref{thm-concentration-vb:expect:regression} and of Corollary~\ref{cor-matrix}, are also provided in the supplementary material.
\end{itemize}

\newpage

\section*{Supplementary material}

\subsection{Consistency of variational approximations in a model where the MAP and the MLE are not defined}

In this subsection, we provide the toy example announced in the paper, where the MLE (and thus the MAP) are not defined. Then, we show that there is a variational approximation that leads to a consistent estimator. This also implies that the tempered posterior is also consistent in this case.

\subsubsection{Description of the model}

We define the statistical model
$$ P_{\theta} = \frac{1}{2}\mathcal{N}(m,\sigma^2) + \frac{1}{2}\mathcal{N}(0,1) $$
for $\theta=(m,\sigma^2)\in\Theta=\mathbb{R} \times (0,1)$. Let $g_\theta$ denote the density of $P_\theta$ with respect to the Lebesgue measure.

The prior $\pi$ is given by: $m\sim\mathcal{N}(0,1)$ and $\sigma^2\sim\mathcal{U}(0,1)$ (uniform distribution).

\subsubsection{Non-existence of the MLE}

It is easy to check that when $X_1,\dots,X_n$ are i.i.d from $P_{(m_0,\sigma_0^2)}$ then the likelihood function
\begin{align*}
L_n(m,\sigma^2)
& = \prod_{i=1}^n g_{(m,\sigma^2)} (X_i)
\\
& = \frac{1}{(2\sqrt{2\pi})^n} \prod_{i=1}^n \left[ \frac{\exp\left(-\frac{(X_i-m)^2}{2\sigma^2}\right)}{\sqrt{\sigma^2}} + \exp\left(-\frac{X_i^2}{2}\right)\right]
\\
& \geq \frac{1}{(2\sqrt{2\pi})^n} \frac{\exp\left(-\frac{(X_i-m)^2}{2\sigma^2}\right)}{\sqrt{\sigma^2}} \exp\left(-\sum_{j\neq i} \frac{X_j^2}{2} \right)
\end{align*}
satisfies, for any $i$,
$$ L_n(X_i,\sigma^2) \xrightarrow[\sigma^2 \rightarrow 0]{} \infty . $$
Thus, the MLE is not defined. For the same reason, the MAP does not exist either.

\subsubsection{A variational approximation family}

We define a family $\mathcal{F}$ that will lead to a consistent estimation. Note that in this case, the variational approximation is not meant to be helpful for computational purposes. On the other hand, it is a very natural one. Indeed, the problem with the MLE is that the likelihood has huge variations around each $X_i$. We will here smooth these variations, creating a kind of maximum of the local mean value of the likelihood. Take the set $\mathcal{F}$ as the set of uniform distributions $\rho_{a,b,c,d}$ over $[a-c,a+c]\times [b-d,b]$, with $(a,b,c,d)\in \mathcal{P} = \{(a,b,c,d)\in\mathbb{R}^4: b\in(0,1), c>0,0<d<b\}$. So we remind that the estimator is then defined as $ \tilde{\pi}_{n,\alpha}({\rm d}\theta|X_1^n) = \rho_{\hat{a},\hat{b},\hat{c},\hat{d} } $ where
\begin{multline*}
(\hat{a},\hat{b},\hat{c},\hat{d})
= \underset{(a,b,c,d)\in\mathcal{P}}{\arg\min} \, \Biggl\{
 - \alpha \int \sum_{i=1}^n \log g_{(m,\sigma^2)}(X_i) \rho_{a,b,c,d}({\rm d}(m,\sigma^2))
 \\
 + \mathcal{K}(\rho_{a,b,c,d},\pi)
 \Biggr\}.
\end{multline*}
Note that the family $\mathcal{F}$ is inspired by Catoni's point of view~\cite{catoni2004statistical} to use a ``perturbed MLE'' in PAC-Bayesian bounds. An application of Theorem~\ref{thm-concentration-vb:expect} leads to the following result.
\begin{prop}
 \label{prop-toy}
 For any $\alpha\in(0,1)$,
 $$ \mathbb{E}  \left[ \int D_{\alpha}(P_{\theta},P_{\theta_0}) \tilde{\pi}_{n,\alpha}({\rm d}\theta|X_1^n) \right]
  \leq 
 \frac{   1.5 \log\left(\frac{2n}{ \sigma_0^2 }\right) + m_0^2 + 1.23}{n(1-\alpha)} . $$
\end{prop}
As a corollary, we also have that the tempered posterior $\pi_{n,\alpha}(\cdot|X_1^n)$ satisfies the same inequality.

\subsubsection{Proof of Proposition~\ref{prop-toy}}

Theorem~\ref{thm-concentration-vb:expect} gives
\begin{multline}
\label{appli-du-thm}
\mathbb{E} \left[ \int D_{\alpha}(P_{\theta},P_{\theta_0}) \tilde{\pi}_{n,\alpha}({\rm d}\theta|X_1^n) \right]
\\
\leq \inf_{(a,b,c,d)\in\mathcal{P}} \Biggl\{ \frac{\alpha}{1-\alpha} \int \mathcal{K}(P_{(m_0,\sigma_0^2)},P_{(m,\sigma)}) \rho_{a,b,c,d}({\rm d}(m,\sigma^2))
 \\
 + \frac{\mathcal{K}(\rho_{a,b,c,d},\pi) }{n(1-\alpha)}
 \Biggr\}.
\end{multline}

We then have:
\begin{align*}
\mathcal{K}(P_{(m_0,\sigma_0^2)},P_{(m,\sigma^2)}) 
& = \mathcal{K}\left(\frac{1}{2}\mathcal{N}(m_0,\sigma_0^2)+\frac{1}{2}\mathcal{N}(0,1),\frac{1}{2}\mathcal{N}(m,\sigma^2)+\frac{1}{2}\mathcal{N}(0,1)\right)
\\
& \leq \frac{1}{2} \mathcal{K}\left(\mathcal{N}(m_0,\sigma_0^2),\mathcal{N}(m,\sigma^2)\right)\text{ (Theorem 11 in~\cite{van2014renyi})}
\\
& = \frac{1}{4}\left[\frac{(m_0-m)^2}{\sigma^2} + \log\left(\frac{\sigma^2}{\sigma_0^2}\right) + \frac{\sigma_0^2-\sigma^2}{\sigma^2} \right].
\end{align*}
Assume that $m\in[m_0-\sqrt{\delta \sigma_0^2},m_0+\sqrt{\delta\sigma_0^2}]$ and $\sigma^2 \in [\sigma^2_0-\delta\sigma^2_0,\sigma^2_0]$ for some $0<\delta<1$. Then:
$$
\mathcal{K}(P_{(m_0,\sigma_0^2)},P_{(m,\sigma^2)})
\leq \frac{1}{4} \left[  \frac{\delta \sigma_0^2}{\sigma^2} + \frac{\delta\sigma_0^2}{\sigma^2} \right] \leq \frac{\delta \sigma_0^2 }{2(\sigma_0^2 - \delta \sigma_0^2)} = \frac{\delta}{2(1-\delta)}.
$$
This implies that for any $\delta\in(0,1)$,
$$
 \int \mathcal{K}(P_{(m_0,\sigma_0^2)},P_{(m,\sigma^2)}) \rho_{m_0,\sigma_0^2,\sqrt{\delta\sigma_0^2},\delta \sigma_0^2}({\rm d}(m,\sigma^2)) \leq \frac{\delta}{2(1-\delta)}.
$$
On the other hand,
\begin{align*}
 \mathcal{K}(\rho_{m_0,\sigma_0^2,\sqrt{\delta\sigma_0^2},\delta\sigma_0^2},\pi)
 & = \mathcal{K}(\mathcal{U}(m_0-\sqrt{\delta \sigma^2},m_0+\sqrt{\delta \sigma^2}),\mathcal{N}(0,1))
 \\
 & \quad \quad + \mathcal{K}(\mathcal{U}(\sigma_0^2-\delta\sigma_0^2,\sigma_0^2),\mathcal{U}(0,1))
 \\
 & = \frac{1}{2\sqrt{\delta \sigma_0^2 }} \int_{m_0-\sqrt{\delta \sigma_0^2}}^{m_0+\sqrt{\delta \sigma_0^2}}  \log\left(\frac{\sqrt{2\pi}}{2\sqrt{\delta\sigma_0^2} \exp\left(\frac{-x^2}{2}\right)} \right) {\rm d}x
 \\ & \quad \quad  + \log\left(\frac{1}{\delta \sigma_0^2} \right)
 \\
 & \leq \frac{1}{2\sqrt{\delta \sigma_0^2 }} \int_{-\sqrt{\delta \sigma_0^2}}^{\sqrt{\delta \sigma_0^2}} \left[ \frac{(m_0+x)^2}{2}+ \log\left(\sqrt{\frac{\pi}{2\delta \sigma_0^2}} \right) \right] {\rm d}x
 \\ & \quad \quad  + \log\left(\frac{1}{\delta \sigma_0^2} \right)
 \\
 & \leq m_0^2 + \delta \sigma_0^2 + \frac{1}{2}\log\left(\frac \pi 2 \right) +  \frac{3}{2} \log\left(\frac{1}{\delta \sigma_0^2 }\right).
\end{align*}
Plugging everything into~\eqref{appli-du-thm} gives:
\begin{multline*}
\mathbb{E} \left[ \int D_{\alpha}(P_{\theta},P_{\theta_0}) \tilde{\pi}_{n,\alpha}({\rm d}\theta|X_1^n) \right]
\\
\leq \inf_{\delta>0} \left[ \frac{\alpha \delta}{2 (1-\alpha)(1-\delta)} +
 \frac{ m_0^2 + \delta \sigma_0^2 + \frac{1}{2}\log\left(\frac \pi 2 \right) +  \frac{3}{2} \log\left(\frac{1}{\delta \sigma_0^2 }\right)}{n(1-\alpha)} \right].
\end{multline*}
The value $\delta = 1/(2n)$ gives, using $(1-\delta)>1/2$ to simplify things,
\begin{multline*}
\mathbb{E} \left[ \int D_{\alpha}(P_{\theta},P_{\theta_0}) \tilde{\pi}_{n,\alpha}({\rm d}\theta|X_1^n) \right]
\\
\leq  \frac{\alpha }{2n (1-\alpha)} +
 \frac{ m_0^2 + \frac{\sigma_0^2}{2n} + \frac{1}{2}\log\left(\frac \pi 2 \right) +  \frac{3}{2} \log\left(\frac{2n}{ \sigma_0^2 }\right)}{n(1-\alpha)}
\\
\leq \frac{0.5 }{n (1-\alpha)} +
 \frac{ m_0^2 + 0.5 + 0.23 +  1.5 \log\left(\frac{2n}{ \sigma_0^2 }\right)}{n(1-\alpha)}.
\end{multline*}

\subsection{Proof of Theorem~\ref{thm-matrix}}

Fix $B>0$, $r\geq 1$ and any pair $(\bar{U},\bar{V})\in\mathcal{M}_{r,B}$ and define for $\delta\in(0,B)$ that will be chosen later,
\begin{equation*}
 \rho_n({\rm d}U,{\rm d}V,{\rm d}\gamma) \propto \mathbf{1}(\|U-\bar{U}\|_{\infty} \leq \delta,\|U-\bar{U}\|_{\infty} \leq \delta) \pi({\rm d}U,{\rm d}V,{\rm d}\gamma).
\end{equation*}
Note that it can be factorized so it belongs to the family $\mathcal{F}$.

We adapt the calculations from~\cite{alquier2014bayesian,alquier2016oracle} but simplify a lot. First, note that
\begin{equation*}
 \mathcal{K}(P_M,P_{UV^t}) = \frac{1}{mp} \sum_{i=1}^m \sum_{j=1}^p \frac{(M_{i,j}-(UV^t)_{i,j})^2}{2\sigma^2}= \frac{\|M-UV^t\|_F^2}{2\sigma^2 mp}
\end{equation*}
and that for any $(U,V)$ in the support of $\rho_n$ we have
\begin{align*}
 \|M-UV^t\|_F & = \|M-\bar{U}\bar{V}^t+\bar{U}\bar{V}^t-\bar{U}V^t+\bar{U}V^t-UV^t\|_F \\
 & \leq  \|M-\bar{U}\bar{V}^t\|_F+\|\bar{U}\bar{V}^t-\bar{U}V^t\|_F+\|\bar{U}V^t-UV^t\|_F \\
 & =  \|M-\bar{U}\bar{V}^t\|_F+\|\bar{U}(\bar{V}^t-V^t)\|_F+\|(\bar{U}-U)V^t\|_F \\
 & \leq \|M-\bar{U}\bar{V}^t\|_F+\|\bar{U}\|_F \|\bar{V}-V\|_F+\|\bar{U}-U\|_F \|V^t\|_F \\
 & \leq  \|M-\bar{U}\bar{V}^t\|_F+mp \|\bar{U}\|_{\infty}^{1/2} \|\bar{V}-V\|_\infty^{1/2}
                       +mp \|V\|_{\infty}^{1/2} \|\bar{U}-U\|_\infty^{1/2} \\
 & \leq \|M-\bar{U}\bar{V}^t\|_F + mp (B^{1/2} \delta^{1/2} + (B+\delta)^{1/2} \delta^{1/2} ) \\
 & \leq \|M-\bar{U}\bar{V}^t\|_F + 2 m p \delta^{1/2} (B+\delta)^{1/2} \\
 & \leq \|M-\bar{U}\bar{V}^t\|_F + 2^{3/2} m p \delta^{1/2} B^{1/2}  = \|M-\bar{U}\bar{V}^t\|_F + B/n
\end{align*}
with the choice $\delta=B/[8(nmp)^2]$ which satisfies $0<\delta<B$.
Then, we derive
\begin{align*}
 \mathcal{K}(\rho_n,\pi) & = \log \frac{1}{\pi\left(\|U-\bar{U}\|_{\infty} \leq \delta,\|U-\bar{U}\|_{\infty} \leq \delta \right)}
 \\
 & =
 \log \frac{1}{\int \pi\left(\|U-\bar{U}\|_{\infty} \leq \delta,\|U-\bar{U}\|_{\infty} \leq \delta|\gamma \right) \pi({\rm d}\gamma) }
 \\
 & =
 \log \frac{1}{\int \pi\left(\|U-\bar{U}\|_{\infty} \leq \delta|\gamma \right) \pi({\rm d}\gamma) } +  \log \frac{1}{\int \pi\left(\|V-\bar{V}\|_{\infty} \leq \delta|\gamma \right) \pi({\rm d}\gamma) }
 \\
 & =  \log \frac{1}{\int_E \pi\left(\|U-\bar{U}\|_{\infty}  \leq \delta|\gamma \right) \pi({\rm d}\gamma) } +  \log \frac{1}{\int_E \pi\left(\|V-\bar{V}\|_{\infty}  \leq \delta|\gamma \right) \pi({\rm d}\gamma) }
\end{align*}
for any event $E$. We actually take $E=\{\gamma_1,\dots,\gamma_r \in[B^2,2B^2],\gamma_{r+1},\dots,\gamma_K \in[s,2s]\}$ and $s\in(0,B^2)$ is to be chosen later.
Then note that
\begin{multline*}
 \pi\left(\|U-\bar{U}\|_{\infty}  \leq \delta|\gamma \right)
  = \pi\left(   \forall i,k: |U_{i,k} - \bar{U}_{i,k}|  \leq \delta|\gamma \right)
 \\
  = \prod_{i=1}^m \prod_{k=1}^r \pi\left(  |U_{i,k} - \bar{U}_{i,k}|  \leq \delta|\gamma_k \right)
      \pi\left( \left. \max_{i=1}^m \max_{k=r+1}^K |U_{i,k}|  \leq \delta\right\vert\gamma_{r+1},\dots,\gamma_K \right).
\end{multline*}
First,
\begin{align*}
 & \pi\left( \left.\max_{i=1}^m \max_{k=r+1}^K |U_{i,k}  |  \leq \delta\right\vert \gamma_{r+1},\dots,\gamma_K \right)
 \\
 & = 1 - \pi\left( \left.\max_{i=1}^m \max_{k=r+1}^K |U_{i,k}|  > \delta\right\vert\gamma_{r+1},\dots,\gamma_K \right)
 \\
 & \geq 1 - \pi\left( \left.\sum_{i=1}^m \sum_{k=r+1}^K |U_{i,k}|  > \delta\right\vert\gamma_{r+1},\dots,\gamma_K \right)
 \\
 & \geq 1 - \frac{\sum_{i=1}^m \sum_{k=r+1}^K \pi\left.\left(|U_{i,k}| \right\vert \gamma_k \right)}{\delta}
 \\
 & \geq 1 - \frac{m (K-r) \max_{k\geq r+1} \sqrt{\gamma_k}}{\delta}
 \\
 & \geq 1 - \frac{m K \sqrt{2 s}}{\delta} = 1/2
\end{align*} with $ s = \frac{1}{2}\left(\frac{\delta}{2(m\vee p)K}\right)^2$
which satisfies $0<s<B^2$. Then, for $k\leq r$,
\begin{align*}
 \pi\left(  |U_{i,k} - \bar{U}_{i,k}|  \leq \delta|\gamma_k \right)
 & = \frac{1}{\sqrt{2\pi\gamma_k}} \int_{\bar{U}_{i,k}-\delta}^{\bar{U}_{i,k}+\delta}
 \exp\left( - \frac{x^2}{2\gamma_k} \right) {\rm d}x
 \\
 & \geq \frac{2\delta \exp\left(-\frac{(B+\delta)^2}{2\gamma_k}\right) }{\sqrt{2\pi\gamma_k}}
 \\
 & \geq  \frac{\delta \exp\left(-\frac{(B+\delta)^2}{2B^2}\right) }{B \sqrt{\pi}} \text{ as } B^2\leq \gamma_k \leq 2 B^2
 \\
 & \geq  \frac{\delta \exp\left(-2\right) }{B\sqrt{\pi}} \text{ as } \delta<1\leq B
\end{align*}
and so
\begin{equation*}
 \prod_{i=1}^m \prod_{k=1}^r \pi\left(  |U_{i,k} - \bar{U}_{i,k}|  \leq \delta|\gamma_k \right)
 \geq  \left( \frac{\delta}{B \sqrt{\pi}}\right)^{mr} \exp\left(-2mr\right) .
\end{equation*}
Finally
\begin{align*}
 \int_E \pi\left(\|U-\bar{U}\|_{\infty}  \leq \delta|\gamma \right) \pi({\rm d}\gamma)
 & \geq  \int_E \frac{1}{2}  \left( \frac{\delta}{B \sqrt{\pi}}\right)^{mr} \exp\left(-2mr\right) \pi({\rm d}\gamma)
 \\
 & = \frac{1}{2} \left( \frac{\delta}{B\sqrt{\pi}}\right)^{mr} \exp\left(-2mr\right) \pi(\gamma\in E)
\end{align*}
and it remains to lower bound
\begin{align*}
 \pi(\gamma\in E)
 & = \left(\prod_{k=1}^{r} \pi(1\leq \gamma_k\leq 2 ) \right) \left(\prod_{k=r+1}^{K} \pi(s\leq \gamma_k\leq 2s ) \right)
 \\
 & = \left( \frac{b^a}{\Gamma(a)} \right)^K \left[ \int_{B^2}^{2B^2}  x^{-a-1}\exp\left(-\frac{b}{x}\right) {\rm d}x \right]^r \left[ \int_{s}^{2s} x^{-a-1}\exp\left(-\frac{b}{x}\right) {\rm d}x\right]^{K-r}
 \\
 & \geq \left( \frac{b^a}{\Gamma(a)} \right)^K  \left[B^2 (2B^2)^{-a-1}\exp\left(-\frac{b}{B^2}\right)\right]^r \left[s(2s)^{-a-1}\exp\left(-\frac{b}{s}\right)\right]^{K-r}
 \\
 & = \left( \frac{b^a}{2^{a+1}\Gamma(a)} \right)^K \exp\left[-\frac{b}{B^2} r - \frac{b}{s}(K-r) \right] (B^{2})^{-(a+1)r} s^{-a(K-r)}
 \\
 & \geq \left( \frac{b^a}{(B^{2})^a 2^{a+1}\Gamma(a)} \right)^K \exp\left[- \frac{Kb}{s}\right]
\end{align*}
as $s<B^2$. So, finally,
\begin{equation*}
  \mathcal{K}(\rho_n,\pi)
   \leq r(m+p) \log\left(\frac{B \sqrt{\pi} \exp(2)}{\delta}\right)
   + K \left[\log\left( \frac{2^{a+1}\Gamma(a)(B^2)^a}{b^a}\right) + \frac{b}{s} \right]+ 2 \log(2).
\end{equation*}
The choice $b=s$ leads to
\begin{align*}
  \mathcal{K}(\rho_n,\pi)
  & \leq  r(m+p) \log\left(\frac{B \sqrt{\pi} \exp(2)}{\delta}\right)
   + K \log\left( \frac{{\rm e} 2^{a+1}\Gamma(a)(B^2)^a}{s^a}\right)+ 2 \log(2)
  \\
  & \leq  r(m+p) \log\left(8 \sqrt{\pi} \exp(2)(nmp)^2\right)
    +4aK \log\left(nmp\right)
    \\ & \quad \quad + K\log({\rm e}2^{10a+1}\Gamma(a))+ 2 \log(2)
\end{align*}
where we replaced $\delta$ and $s$ by their respective value. In order to keep the expressions as simple as possible we can use $K\leq m\vee p \leq m+p \leq r(m+p)$ and $2\leq m+p\leq r(m+p)$ to get
\begin{equation*}
  \mathcal{K}(\rho_n,\pi) \leq 2(1+2a) r(m+p) \left[ \log(nmp)
  +
  \underbrace{\log(8\sqrt{\pi}\Gamma(a)2^{10a+1})+3}_{=: \mathcal{C}(a)}
  \right].
\end{equation*}

We are now in position to apply Theorem \ref{thm-concentration-vb:expect}. Then
\begin{multline*}
\mathbb{E} \left[ \int D_{\alpha}(P_{\theta},P_{\theta_0}) \tilde{\pi}_{n,\alpha}({\rm d}\theta|X_1^n) \right]
\\
\leq  \frac{\alpha}{1-\alpha} \int \mathcal{K}(P_{\theta_0},P_{\theta}) \rho_n({\rm d}\theta)
 + \frac{\mathcal{K}(\rho_n,\pi) }{n(1-\alpha)}
\\
\leq \frac{\alpha}{1-\alpha} \frac{\left[\|M-\bar{U}\bar{V}^t\|_F + \frac{\sqrt{B}}{n}\right]^2}{2\sigma^2 mp}
 \\
 + \frac{2(1+\alpha)(1+2a) r(m+p) \left[ \log(nmp)
  + \mathcal{C}(a,B)
  \right]}{n(1-\alpha)}.
\end{multline*}

\subsection{Proof of Corollary~\ref{cor-matrix}}

We start from~\eqref{pre-coro-matrix}. Under the boundedness assumption on $M_0$ it is obvious that $\forall M$, $d_{\alpha,\sigma}(M,M_0) \geq d_{\alpha,\sigma}({\rm clip}_C(M),M_0)$,
so
\begin{multline}
\label{frejods}
\mathbb{E} \left[ \int d_{\alpha,\sigma}({\rm clip}_C(M),M_0) \tilde{\pi}_{n,\alpha}({\rm d}M|X_1^n) \right]
 \\
 \leq
  \frac{2(1+\alpha)(1+2a) r(m+p) \left[ \log(nmp)
  + \mathcal{C}(a) + \frac{\alpha B}{2\sigma^2 mp}
  \right]}{n(1-\alpha)}.
\end{multline}
Fix $M$ and for short, put $N={\rm clip}_C(M)$. We have:
\begin{align*}
 d_{\alpha,\sigma}(N,M_0)
 & =
 \frac{-1}{1-\alpha} \log \left[ \frac{1}{mp}\sum_{i=1}^m\sum_{j=1}^p \exp\left(\frac{\alpha(\alpha-1)(N_{i,j}-(M_0)_{i,j})^2}{2\sigma^2}\right) \right]
 \\
 & \geq
 \frac{1}{1-\alpha} \left[ 1 - \frac{1}{mp}\sum_{i=1}^m\sum_{j=1}^p \exp\left(\frac{\alpha(\alpha-1)(N_{i,j}-(M_0)_{i,j})^2}{2\sigma^2}\right) \right].
\end{align*}
By assumption, $(N_{i,j}-(M_0)_{i,j})^2 / (2 \sigma^2) \leq (2C)^2 / (2\sigma^2) = 2 C^2 / \sigma^2 $. Straightforward derivations show that for any $x\in[0,2 C^2 / \sigma^2]$ we have
$$ 1 - \left(  \frac{\sigma^2[1-\exp(2 C^2 \alpha(\alpha-1) / \sigma^2)]}{2 C^2}\right) x \geq \exp(\alpha(\alpha-1)x), $$
this leads to
$$
 d_{\alpha,\sigma}(N,M_0)
\geq
 \frac{1}{1-\alpha}  \frac{\sigma^2[1-\exp(2 C^2 \alpha(\alpha-1) / \sigma^2)]}{2 C^2} \frac{\|N-M_0\|_F^2}{2\sigma^2}.
$$
Pluging this into~\eqref{frejods} gives the result claimed.

\subsection{Proof of Theorem~\ref{thm-concentration-vb:expect:regression}}

Let $(\beta_k^0)$ denote the coefficients of $f^0$: $f_0 = \sum_{k=1}^\infty \beta_k^0 \varphi_k$. Theorem~\ref{thm-concentration-vb:expect} gives:
\begin{align*}
\mathbb{E} & \left[ \int D_{\alpha}(P_{f},P_{f_0}) \tilde{\pi}_{n,\alpha}({\rm d}\theta|X_1^n) \right]
\\
& \leq \frac{1}{1-\alpha} \inf_{K\geq 1} \inf_{\mathbf{m},s^2} \Biggl\{ \frac{\alpha}{2} \int \int_{-1}^1 \left(f_0 - \sum_{k=1}^K \beta_k \varphi_{k} \right)^2 \Phi({\rm d}\beta,\mathbf{m},s^2)
 \\
 & \quad \quad + \frac{\sum_{k=1}^K \frac{1}{2}  \left[\log\left(\frac{1}{s^2}\right)+s^2+m_k^2-1\right]  + K\log(2)}{n} \Biggr\}
\\
& = \frac{1}{1-\alpha} \inf_{K\geq 1} \inf_{\mathbf{m},s^2} \Biggl\{ \frac{\alpha}{2} \int \int_{-1}^1 \left(\sum_{k=1}^\infty \beta_k^0 - \sum_{k=1}^K \beta_k \varphi_{k} \right)^2 \Phi({\rm d}\beta,\mathbf{m},s^2)
 \\
 & \quad \quad + \frac{\sum_{k=1}^K \frac{1}{2}  \left[\log\left(\frac{1}{s^2}\right)+s^2+m_k^2-1\right]  + K\log(2)}{n} \Biggr\}
\\
& = \frac{1}{1-\alpha} \inf_{K\geq 1} \inf_{\mathbf{m},s^2} \Biggl\{ \frac{\alpha}{2} \sum_{k=1}^K (m_k-\beta_k^0)^2 + \frac{\alpha K s^2}{2} + \frac{\alpha}{2} \sum_{k=K+1}^\infty (\beta_k^0)^2
 \\
 & \quad \quad + \frac{\sum_{k=1}^K \frac{1}{2}  \left[\log\left(\frac{1}{s^2}\right)+s^2+m_k^2-1\right]  + K\log(2)}{n} \Biggr\}.
\end{align*}
The choice $(m_1,\dots,m_K)=(\beta_1,\dots,\beta_K)$ gives:
\begin{multline*}
\mathbb{E}  \left[ \int D_{\alpha}(P_{f},P_{f_0}) \tilde{\pi}_{n,\alpha}({\rm d}\theta|X_1^n) \right]
 \leq \frac{1}{1-\alpha} \inf_{K\geq 1}  \Biggl\{  \frac{\alpha K s^2}{2} + \frac{\alpha}{2} \sum_{k=K+1}^\infty (\beta_k^0)^2
 \\
  + \frac{ \sum_{k=1}^K \frac{1}{2} \left[\log\left(\frac{1}{s^2}\right)+s^2+(\beta_k^0)^2-1\right]  + K\log(2)}{n} \Biggr\}.
\end{multline*}
From Chapter 1 in~\cite{tsybakov2009book} we know that $f_0 \in\mathcal{W}(r,C^2)$ implies
$ \sum_{k=K+1}^\infty (\beta_k^0)^2 \leq \Lambda(r,C) K^{-2r} $ for some $\Lambda(k,C)$. Moreover, $\sum_{k=1}^{K}  (\beta_k^0)^2 \leq \sum_{k=1}^\infty k^{2r} (\beta_k^0)^2 \leq C^2 $. So finally:
\begin{multline*}
\mathbb{E}  \left[ \int D_{\alpha}(P_{f},P_{f_0}) \tilde{\pi}_{n,\alpha}({\rm d}\theta|X_1^n) \right]
 \leq \frac{1}{1-\alpha} \inf_{K\geq 1} \Biggl\{ \frac{\alpha \Lambda(r,C)}{2} K^{-2r}
 \\
 + \frac{\alpha K s^2}{2} + \frac{ K\left[ \log(2) + \frac{s^2}{2} + \frac{\log(\frac{1}{s^2})}{2} \right] +C^2 }{n} \Biggr\}.
\end{multline*}
The choice $s^2 = 1/n$ leads to
\begin{multline*}
\mathbb{E}  \left[ \int D_{\alpha}(P_{f},P_{f_0}) \tilde{\pi}_{n,\alpha}({\rm d}\theta|X_1^n) \right]
 \leq \frac{1}{1-\alpha} \inf_{K\geq 1} \Biggl\{ \frac{\alpha \Lambda(r,C)}{2} K^{-2r}
 \\
  + \frac{ K\left[ \log(2) + \frac{\alpha}{2} + \frac{1}{2n} + \frac{\log(n)}{2} \right] +C^2 }{n} \Biggr\}.
\end{multline*}
The choice $K = \lceil (n/\log(n))^{1/(2r+1)} ) \rceil $ leads to the result.

\end{document}